\documentclass[11pt, a4paper, leqno]{article}
\usepackage[utf8]{inputenc}
\usepackage{amsmath}
\usepackage{amsfonts, amssymb}
\usepackage{scalerel}
\usepackage{afterpage, pdflscape} 
\usepackage{booktabs, multirow, array, caption} 
\usepackage{tikz} 
\definecolor{LightCyan}{rgb}{0.88,1,1}
\definecolor{Palla}{RGB}{255,255,189}
\usepackage{xcolor}
\usepackage{subfigure}
\usepackage{mathtools}
\usepackage{graphicx}
\usepackage{epstopdf}
\usepackage{bm} 
\usepackage{xspace} 
\usepackage[sort&compress, numbers]{natbib}
\usepackage{fullpage}
\usepackage{authblk}
\usepackage{titlesec}
\usepackage[titletoc,toc,title]{appendix}
\usepackage{rotating,lineno}
\usepackage{caption}

\newcommand\bbox{\scaleobj{0.6}{\Box}}

\newcommand\ys{$y^{\bbox}_1$\xspace}
\newcommand\yg{$y_1$\xspace}
\newcommand\ygg{$y_2$\xspace}
\newcommand\zs{$z^{\bbox}$\xspace}
\newcommand\zg{$z$\xspace}

\newcommand\mys{y^{\bbox}_1}
\newcommand\myg{y_1}
\newcommand\mygg{y_2}
\newcommand\mzs{z^{\bbox}}
\newcommand\mzg{z}

\newcommand{\alps}{$\alpha^{\bbox}$\xspace}
\newcommand{\alpg}{$\alpha$\xspace}
\newcommand{\gammas}{$\gamma^{\bbox}_1$\xspace}
\newcommand{\gammag}{$\gamma_1$\xspace}
\newcommand{\gammagg}{$\gamma_2$\xspace}
\newcommand{\ks}{$\kappa^{\bbox}$\xspace}
\newcommand{\kg}{$\kappa_1$\xspace}
\newcommand{\kgg}{$\kappa_2$\xspace}
\newcommand{\nus}{$\nu^{\bbox}$\xspace}
\newcommand{\nug}{$\nu$\xspace}
\newcommand{\zets}{$\zeta^{\bbox}$\xspace}
\newcommand{\zetg}{$\zeta$\xspace}

\newcommand{\malps}{\alpha^{\bbox}}
\newcommand{\malpg}{\alpha}
\newcommand{\mgammas}{\gamma^{\bbox}_1}
\newcommand{\mgammag}{\gamma_1}
\newcommand{\mgammagg}{\gamma_2}
\newcommand{\mks}{\kappa^{\bbox}}
\newcommand{\mkg}{\kappa_1}
\newcommand{\mkgg}{\kappa_2}
\newcommand{\mnus}{\nu^{\bbox}}
\newcommand{\mnug}{\nu}
\newcommand{\mzets}{\zeta^{\bbox}}
\newcommand{\mzetg}{\zeta}

\makeatletter
\newcommand{\pushright}[1]{\ifmeasuring@#1\else\omit\hfill$\displaystyle#1$\fi\ignorespaces}
\newcommand{\pushleft}[1]{\ifmeasuring@#1\else\omit$\displaystyle#1$\hfill\fi\ignorespaces}
\makeatother

\title{Host-virus evolutionary dynamics with specialist and generalist infection strategies: bifurcations, bistability and chaos}

\author[1,2]{Anel~Nurtay}
\author[1,3]{Matthew~G.~Hennessy}
\author[4,1,5]{Llu\'is~Alsed\`a}
\author[6,7,8]{Santiago~F.~Elena}
\author[1,5,8]{Josep~Sardany{\'e}s}

\affil[1]{Centre de Recerca Matem\`{a}tica, Campus de Bellaterra, Edifici C, 08193 Bellaterra, Spain}
\affil[2]{Big Data Institute, Univeristy of Oxford, Old Road Campus, Oxford, OX3 7LF, United Kingdom}
\affil[3]{Mathematical Institute, Univeristy of Oxford, Andrew Wiles Building, Radcliffe Observatory Quarter, Woodstock Road, Oxford, OX2 6GG, United Kingdom}
\affil[4]{Departament de Matem\`{a}tiques, Universitat Aut\`{o}noma de Barcelona, Campus de Bellaterra, Edifici C, 08193 Bellaterra, Spain}
\affil[5]{Barcelona Graduate School of Mathematics (BGSMath). Campus de Bellaterra, Edifici C, 08193 Bellaterra, Spain}
\affil[6]{Instituto de Biolog\'ia Integrativa de Sistemas, CSIC-Universitat de Val\`encia, Parc Cient\'ific UV, Paterna, 46182 Val\`encia, Spain}
\affil[7]{The Santa Fe Institute, 1399 Hyde Park Road, Santa Fe NM87501, USA}
\affil[8]{Dynamical Systems and Computational Virology, Associated Unit Instituto de Biolog\'ia Integrativa de Sistemas-Centre de Recerca Matem\`{a}tica}

\newcommand*\patchAmsMathEnvironmentForLineno[1]{%
 \expandafter\let\csname old#1\expandafter\endcsname\csname #1\endcsname
 \expandafter\let\csname oldend#1\expandafter\endcsname\csname end#1\endcsname
 \renewenvironment{#1}%
    {\linenomath\csname old#1\endcsname}%
    {\csname oldend#1\endcsname\endlinenomath}}%
\newcommand*\patchBothAmsMathEnvironmentsForLineno[1]{%
 \patchAmsMathEnvironmentForLineno{#1}%
 \patchAmsMathEnvironmentForLineno{#1*}}%
\AtBeginDocument{%
\patchBothAmsMathEnvironmentsForLineno{equation}%
\patchBothAmsMathEnvironmentsForLineno{align}%
\patchBothAmsMathEnvironmentsForLineno{flalign}%
\patchBothAmsMathEnvironmentsForLineno{alignat}%
\patchBothAmsMathEnvironmentsForLineno{gather}%
\patchBothAmsMathEnvironmentsForLineno{multline}%
}

\providecommand{\keywords}[1]{\textbf{Keywords: } #1}

\begin{document}
	
\maketitle

\begin{abstract}
In this work we have investigated the evolutionary dynamics of a generalist pathogen, {\it e.g.} a virus population, that evolves towards specialisation in an environment with multiple host types.  We have particularly explored under which conditions generalist viral strains may rise in frequency and coexist with specialist strains or even dominate the population.  By means of a nonlinear mathematical model and bifurcation analysis, we have determined the theoretical conditions for stability of nine identified equilibria and provided biological interpretation in terms of the infection rates for the viral specialist and generalist strains.  By means of a stability diagram we identified stable fixed points and stable periodic orbits, as well as regions of bistability.  For arbitrary biologically feasible initial population sizes, the probability of evolving towards stable solutions is obtained for each point of the analyzed parameter space.  This probability map shows combinations of infection rates of the generalist and specialist strains that might lead to equal chances for each type becoming the dominant strategy.  Furthermore, we have identified infection rates for which the model predicts the onset of chaotic dynamics. 
Several degenerate Bogdanov-Takens and zero-Hopf bifurcations are detected along with generalized Hopf and zero-Hopf bifurcations. This manuscript provides additional insights into the dynamical complexity of host-pathogen evolution towards different infection strategies.
\end{abstract}

\keywords{mathematical modelling; within-host dynamics; generalist and specialist; bifurcation analysis, stability}


\section{Introduction}
Within-host evolution of microparasites can burden the treatment of a disease and lead to the rise of drug-resistant strains or the emergence of new strains with enhanced epidemic potential. Polymorphism of microparasites within an individual host is not limited to cases of rapidly evolving RNA viruses~\cite{sanjuan2012molecular, jenkins2002rates} such as influenza A~\cite{leonard2016deep,hom2019deep}  or HIV-1~\cite{fraser2014virulence, lemey2006hiv} and DNA viruses~\cite{simmonds2018prisoners, duffy2008rates} such as parvovirus~\cite{shackelton2005high, young2004parvovirus}.  Bacterial pathogens, due to their larger genomes, can have similarly high per-genome evolutionary rates~\cite{didelot2016within}. The changes occurring in microparasites within a host can contribute to optimization of the overall fitness of the microparasite, and can tailor it to align with the changes in the host responses to infection. However, adaptation towards a new trait in a host or a change of the host type altogether comes with a price for infectious agents~\cite{bedhomme2015emerging}. As an illustration, Dengue virus (DENV) variants --- an arbovirus that circulates between humans and mosquitoes --- are positively or negatively selected in each of their alternative hosts 
leading to drastically differing populations of the virus in each host type~\cite{villordo2015dengue}. Every reintroduction of DENV to human populations threatens to bring variants that carry new characteristics. Due to the sensitivity of the pathogen's diversity on available host resources, the ability to quantify the dependence of pathogen evolution on its diverse environment is highly relevant to understanding the onset of future epidemics.

The ease of interpretation and tractability of Lotka--Volterra-type models~\cite{yasuhiro1996global} has made them a valuable tool in the investigation of interacting populations. Due to similarity, research on host-pathogen dynamics successfully employed methods and results from predator-prey-type models~\cite{roy2008effects, begon1995beyond, anderson1982coevolution, hethcote1989periodicity}. 
A single disease occurring in a population of one host type has been modeled by Anderson and May~\cite{anderson1979population}, providing a setup for various modifications of their model. 
For example, models with two host types and a homogeneous pathogen population were used to study the transmission of a disease between host types~\cite{holt1985infectious, begon1992disease, gandon2004evolution}. Conceptually opposite problems, that aimed to model coinfection of a host with different types of pathogens, have also been proposed~\cite{gog2002status, bianco2011asymmetry, nyabadza2015implications}. Similar to coinfection models, mathematical models describing within-host dynamics of microparasites~\cite{hoshen2000mathematical, antia1994within}, along with providing immunological insights, have helped to understand the polymorphism of the pathogens' populations.

The specialization of an ancestral generalist pathogen into a new host is the main phenomenon modeled and studied in this work. We use differential equations to study the dynamics of interacting host cells and virus. To model specialization of a viral strain, multiple types of host cells must be presented. Specifically, we introduce two types of susceptible cells as hosts for  an evolving viral population. Part of the virus population is able to infect both susceptible cell types with equal efficiency and is refered to as the generalist strain. Some infections of a designated type of susceptible cells are assumed to benefit the emergence of a specialist strain which, according to the trade-off hypothesis~\cite{bedhomme2015emerging, whitlock1996red}, lacks the ability to infect the second cell type. Under these assumptions, we study the long-term behavior of the system for varying infection rates and present theoretically and numerically obtained conditions that produce qualitatively different dynamical outcomes. 

We have investigated the dependence of the system dynamics on parameters by means of bifurcation analyses~\cite{li2012joint, kuznetsov1994nonlinear}. For nonlinear and parameter-rich models, there is an advantage of conducting parameter analyses as opposed to fixing the whole parameter set~\cite{alizon2010within}. Firstly, because the same dynamics can occur for an interval of parameter values. Secondly, inferring real parameter values for viral strains
can be a resource- and time-demanding goal~\cite{drummond2003inference}. Bifurcation analysis can not only relax requirements for exact parameter values but also determine parameter combinations that lead to certain dynamics of interest. For instance, in a study of the coexistence of wild-type and mutant viral strains \cite{nurtay2019theoretical}, it was shown that the wild-type cannot exist without its mutant, and 
the theoretical lower bound of the mutation rate that leads to a complete elimination of the viral strains from the system was obtained. Another use of bifurcation analysis yielded the characteristics of the virus that, in the absence of an immune response, can eradicate a growing tumor~\cite{wu2004analysis}. In the context of West Nile virus propagation~\cite{blayneh2010backward}, a numerical bifurcation study revealed that the mosquito-reduction strategies can lead to better results than personal human-mosquito interaction control. As a more general case, a model of two parasitoids placed in a two-dimensional environment consisting of two host types has been investigated by means of bifurcation analysis and has yielded heterogeneous spatiotemporal patterns of host and parasitoid abundances~\cite{pearce2006modelling}; these patterns were found to be driven by periodic traveling waves. When phenotype-dependent births and deaths are involved in the dynamics of competing and changing populations, periodic behavior can be one of the likely outcomes~\cite{dieckmann1995evolutionary, dercole2003bifurcation}. Similar to these studies, we here employ bifurcation analysis to describe the dependency of the dynamics on the traits (\emph{i.e.}\ infectiousness) of generalist and specialist strains. We also provide theoretical expressions for bifurcation curves (when applicable) and confirm the two-dimensional bifurcation diagram by means of numerical results that include comparing the bifurcation diagram with an independently obtained probabilistic stability diagram.

The paper is organized as follows. In Section \ref{sec:model}, we describe the mathematical model under study. 
The equilibria of the model and corresponding stability analyses are provided in Section~\ref{sec:equil}.
To further investigate the dynamics, a numerical bifurcation analysis is conducted in Section~\ref{sec:onedim}, revealing regions of bistability and chaotic dynamics. 
The regions of stability (monostability and bistability) are studied using a probabilistic approach and plotted in terms of the infection rates of the generalist and specialist in 
Section~\ref{sec:prob}. The chaotic dynamics are studied in Section~\ref{sec:chaos}. Finally, in Section~\ref{sec:conc}, we summarize and discuss the main results of this work.

\section{Host-pathogen mathematical model with different infection strategies}
\label{sec:model}

In this section we introduce the mathematical model used to investigate the dynamics of cell infection by viruses with generalist and specialist infection strategies. For the sake of clarity, the model is introduced below step by step. Depending on the context, and for simplicity, the variables $x_1, x_2, \mys, \myg, \mygg, \mzs, \mzg$ will be used to denote specific species or their corresponding population sizes.  Figure~\ref{diag} shows a schematic diagram of the dynamical system. 

Consider two different types of cells, $x_1$ and $x_2$, that are placed in a finite environment with carrying capacity $K>0$. The uninfected cells $x_1$ and $x_2$ have growth and death rates $\beta_1 > 0$, $\delta_1 > 0$ and $\beta_2 >0$, $\delta_2 >0$, respectively. Throughout this paper, we assume that $\beta_1 \neq \beta_2$ and $\delta_1 \neq \delta_2$ in order to differentiate the two cell populations. To ensure the persistence of the uninfected cells in the absence of virions, we also assume that $\beta_i > \delta_i$ for $i = 1, 2$. 
Both types of cells can be infected by a strain of virus \zg, hereafter referred to as the generalist, with rate \alpg $>0$.  However, only one of the cell types, $x_1$, can be infected by \zs, the specialist, with rate \alps $> 0$. The rate of infection is considered to be a property of the virus strain, and is different from the incident rate in a population of susceptible organisms. An analogue of \alpg and \alps in a simple epidemiological model with a fixed area of occupation would be the transmission rate~\cite{begon2002clarification}. The dynamics of the cell populations are modeled as follows:
\begin{subequations}
\label{eqn:dim_all}
\begin{align}\label{odes_twocell_spec_gen_x}
    \dot{x}_1 &= \beta_1\,x_1\left(1-\frac{x_1+x_2}{K}\right)-x_1\,\malpg\,\mzg-x_1\,\malps\,\mzs - \delta_1\,x_1,\\
	\dot{x}_2 &= \beta_2\,x_2\left(1-\frac{x_1+x_2}{K}\right)-x_2\,\malpg\,\mzg-\delta_2\,x_2.
\end{align}

\begin{figure}[ht] 
\centering
\captionsetup{width=\linewidth}
\includegraphics[width=0.85\linewidth]{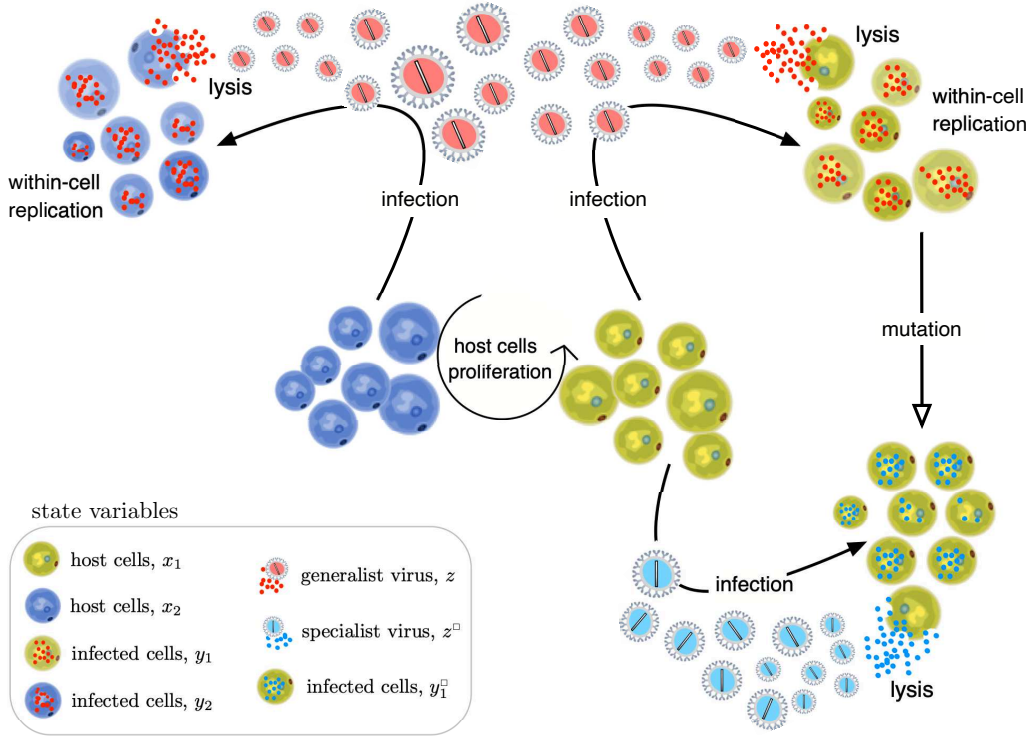}
\caption{Schematic diagram of the studied dynamical system considering two different types of 
host cells, $x_{1,2}$, susceptible to infection.
Two different viral strains are also considered: a generalist strain, $z$, able to infect both types of healthy cells, giving place to infected cells $y_{1,2}$; and a specialist strain, $\protect\mzs$, capable of only infecting cells $x_1$, giving place to cell populations $\protect\mys$. Note that the generalist strain can mutate (open black arrow) to the specialist one. This process is assumed to occur during within-cell virus replication. Recall that throughout this article the variables are also used to label the population sizes (\emph{e.g.} $z$ can refer either to the generalist strain or to its population size).}\label{diag}
\end{figure}

Hereafter, $x_2$ is referred to as a second cell type but can be interpreted as combination of more than one cell type that are not available for the specialist. 
Uninfected cells of type $x_1$ after being infected by \zs become infected cells \ys; uninfected cells of type $x_1$ and $x_2$ after being infected by \zg become infected cells \yg and \ygg, respectively. The growth of the total population of infected cells is assumed to be driven by the infection only, that is, the infected cells do not proliferate. 
The effect of infection on a cell population is depicted by the change in death rates of the infected cell types. Moreover, the effect of infection from a particular virus is also captured and, overall, the death rates of infected cells \ys, \yg, and \ygg are \gammas $>0$, \gammag $>0$, and \gammagg $>0$, respectively.  It is assumed that the generalist strain can specialize towards the infection of only cell type $x_1$.  We thus consider mutation of the generalist \zg into the specialist \zs with rate $\mu>0$.  This is modelled by the transfer from population of \yg to \ys since most of the mutations of viral strains occur during replication and transcription processes inside an infected cell. 
The dynamics of infected cells can be modeled as follows: 
\begin{align}
	\dot{\mys}&=\malps\,\mzs\,x_1 +\mu \myg-\mgammas\,\mys,\\
	\dot{\myg}&=\malpg\,\mzg\,x_1 -\mu \myg-\mgammag\,\myg,\\
	\dot{\mygg}&=\malpg\,\mzg\,x_2 -\mgammagg\,\mygg.
\end{align}
The specialist-infected cells \ys have a burst size \ks $>0$, while generalist-infected cells \yg and \ygg have burst sizes \kg $>0$ and \kgg $>0$, respectively. We let  \nus $>0$ and \nug $>0$ denote the multiplicity of infection (MOI) of \zs and \zg, respectively. Viral populations have a decay rate \zetg $>0$ for the generalist strain and \zets $>0$ for the specialist strain. Under these assumptions, the changes in sizes of the viral populations can be described as follows:
\begin{align}
	\dot{\mzs}&=\mks\mgammas\,\mys-\mnus\,\malps\,\mzs\,x_1-\mzets\,\mzs,\\
	\dot{\mzg}&=\mkg\,\mgammag\,\myg+\mkgg\mgammagg\,\mygg-\mnug\,\malpg\,\mzg(x_1+x_2)-\mzetg\,\mzg.
\end{align}
\end{subequations}
Hereafter, the vector notation for the variables in Eqs.~(1) (as well as in Eqs.~\eqref{odes_twocell_spec_gen_nd}, see below), will be given in the following order: 
$\bm{v} = (x_1,x_2,\mys,\myg,\mygg,\mzs,\mzg)$. The initial population will consist of host cells and virions but no infected cells; that is, $\bm v(0) = (x_1(0)>0, x_2(0)>0, \mys(0) = 0, \myg(0) = 0, \mygg(0) = 0,\mzs(0) > 0,\mzg(0) > 0)$. 

The nondimensionalization of system \eqref{eqn:dim_all} is performed by introducing dimensionless variables that are based on characteristic time scales and population sizes. The quantity $\beta_1 - \delta_1$ describes the effective growth rate of the first uninfected cell type and its inverse, $(\beta_1 - \delta_1)^{-1}$, is used to define the characteristic time scale of the system. 
In a virus-free environment, the total maximum size of the uninfected cell population $x_1$ is  $x_{\text{max}}=(1-\delta_1 / \beta_1)\,K$. Without loss of generality, we choose the characteristic population size to be $x_\text{max}$ for both uninfected and infected cells. Similarly, the characteristic population size for the viruses is chosen to be $x_{\text{max}} \kappa^{\bbox}$, where 
the burst size is required due to a difference in measurement units of the viral loads and the cell populations in the system. The choice to scale with the specialist burst size is arbitrary.
We therefore nondimensionalize the variables according to 
\begin{equation}
t = (\beta_1 - \delta_1)^{-1}\bar t, \quad  x_i = x_\text{max}\bar x_i, \quad y^j_i =  x_\text{max}\,\bar y^j_i,  \quad z^j = x_\text{max}\kappa^{\bbox}\,\bar{z}^{j},
\nonumber\\
\end{equation}
where $i = 1, 2$ and the superscript $j = \bbox$ or its absence represents quantities associated with the speciality and generalist, respectively.
The new variables $\bar t, \bar{x}_i, \bar{y}^j_i, \bar{z}^j$ are non-dimensional versions of time and the population sizes of uninfected cells, infected cells, and virus strains, respectively. 
This rescaling also introduces non-dimensional parameters given by
\begin{align}
\bar \beta_1 = \frac{\beta_2-\delta_2}{\beta_1-\delta_1},
\quad 
\bar \beta_2 = \frac{\beta_2}{\beta_1}, 
\quad 
\bar \alpha = \frac{\kappa^{\bbox}\,K}{\beta_1}{\alpha},
\quad 
\bar \alpha^{\bbox} = \frac{\kappa^{\bbox}\,K}{\beta_1}{\alpha^{\bbox}},
\nonumber \\
\bar \nu = \frac{\nu}{\kappa^{\bbox}}, 
\quad 
\bar{\nu}^{\bbox} = \frac{\mnus}{\kappa^{\bbox}}, 
\quad 
\bar \kappa_i = \frac{\kappa_i}{\kappa^{\bbox}},
\quad 
\bar p = \frac{p}{\beta_1-\delta_1}, \nonumber
\end{align}
where $p$ stands for all other parameters, namely: $\mu$, \gammas, \gammag, \gammagg, \zets, \zetg. Notice that the non-dimensional parameters still remain positive since we generically set $\beta_i>\delta_i, i = 1,2$. 
From a biological point of view, the MOI is several orders of magnitude lower than the burst size. For instance, averaged burst sizes obtained from \emph{in vivo} estimations for simian immunodeficiency virus were about $4 \times 10^4$ and $5.5 \times 10^4$ virions per cell~\cite{chen2007determination}. Concerning the MOI, although experimental estimates are scarce and difficult to obtain, values of about 2-3 have been reported for different DNA or RNA bacteriophages \cite{horiuchi19752,  olkkonen1989quantitation, susskind1974superinfection, turner1999hybrid}, while a MOI value of 4 to 5 was reported for \emph{Autographa californica} nuclear polyhedrosis virus infecting larvae of the moth \emph{Tricoplusia ni} \cite{bull2001persistence}. Also a MOI of about 3 was inferred from the number of proviral copies of HIV-1 in spleen cells of infected patients \cite{jung2002recombination}. Since the burst sizes are extremely large compared to the MOIs, we can assume that $\bar{\nu}, \bar{\nu}^{\bbox} < 1$ and $\bar{\nu}, \bar{\nu}^{\bbox} < \bar{\kappa}_{1,2}$.

Summarizing, by omitting the bars for simplicity, we obtain the following nondimensionalized system of dynamical equations:
\begin{equation}\label{odes_twocell_spec_gen_nd}
\left\{\ 
\begin{split}
	\dot{x}_1 &= x_1\left(1-x_1-x_2\right)-x_1\,\malpg\,\mzg-x_1\,\malps\,\mzs,\\ 
	\dot{x}_2 &= x_2\left[\beta_1-\beta_2(x_1+x_2)\right]-x_2\,\malpg\,\mzg,\\ 
	\dot{\mys}&=\malps\,\mzs\,x_1 +\mu \myg-\mgammas\,\mys,\\
	\dot{\myg}&=\malpg\,\mzg\,x_1 -\mu \myg-\mgammag\,\myg,\\
	\dot{\mygg}&=\malpg\,\mzg\,x_2 -\mgammagg\,\mygg,\\
	\dot{\mzs}&=\mgammas\,\mys-\mnus\,\malps\,\mzs\,x_1-\mzets\,\mzs,\\
	\dot{\mzg}&=\mkg\,\mgammag\,\myg+\mkgg\mgammagg\,\mygg-\mnug\,\malpg\,\mzg(x_1+x_2)-\mzetg\,\mzg.
\end{split}	
\right.	
\end{equation}


\section{Equilibria, linear, and nonlinear stability analysis}
\label{sec:equil}
To gain insight into the general behavior of the model \eqref{odes_twocell_spec_gen_nd}, we calculate equilibrium solutions, conduct a linear stability analysis, and identify bifurcation points.  Due to the high dimensionality of the system, some bifurcation points do not lead to a change in stability because some non-critical eigenvalues retain their positive real part.  There are a number of codimension two bifurcations that occur and these are classified in terms of their normal form coefficients.  MATCONT software \cite{dhooge2003matcont, dhooge2008new} is used to perform numerical continuation of equilibria, periodic solutions, and bifurcation points when required.

\subsection{Equilibrium solutions}

The critical (equilibrium) points of the nondimensionalized system are the solutions 
of the following nonlinear system of equations:
\begin{align}
& \begin{cases}
  x_1\left(1 - x_1 - x_2 - \alpha z - \malps \mzs\right) = 0,\\
  x_2\left(\beta_1 - \beta_2(x_1+x_2) - \alpha z\right) = 0,\\ 
  \kappa_1\gamma_1 y_1 + \kappa_2\gamma_2 y_2 - z \left( \nu\alpha (x_1+x_2) + \zeta\right) = 0,\\
  \gamma^{\bbox}_1 \mys - \mzs \left(\nu^{\bbox} \malps x_1 + \zeta^{\bbox} \right) = 0,
\end{cases} \label{eqn:equil_1}\\
& \begin{cases}  
 \alpha x_1 z - \mu y_1 - \gamma_1 y_1 = 0,\\
 \alpha x_2 z - \gamma_2 y_2 = 0,\\
 \malps x_1 \mzs + \mu y_1 - \gamma^{\bbox}_1 \mys = 0.
\end{cases} \label{eqn:equil_2}
\end{align}
It is immediately clear that the trivial solution 
${\bm{v}}_0 = (0, 0, 0, 0, 0, 0, 0)$ 
is an equilibrium point and it corresponds to a state 
whereby all of the populations become extinct. 
The other equilibria can be found by first solving the three equations
in \eqref{eqn:equil_2} for the populations of the infected cells:
$y_1,\ y_2$ and $\mys$,
and substituting the results into \eqref{eqn:equil_1} to find that:
\begin{gather}
  \hspace*{-10em}
  \left.\begin{cases}
   x_1\left(1 - x_1 - x_2 - \alpha z - \malps \mzs\right) = 0\\
   x_2\left(\beta_1 - \beta_2(x_1+x_2) - \alpha z\right) = 0\\
   z\left(A x_1  + B x_2 - \zeta \right) = 0\\
   C x_1 z - \mzs\left(D x_1 + \zeta^{\bbox}\right) = 0
\end{cases}\mkern-15mu\right\}\quad\Longleftrightarrow \nonumber \\
\hspace*{10em}\begin{cases}
   x_1 = 0\quad\text{or}\quad 1 - x_1 - x_2 - \alpha z - \malps \mzs = 0,\\
   x_2 = 0\quad\text{or}\quad \beta_1 - \beta_2(x_1+x_2) - \alpha z = 0,\\
   \hspace*{0.5em}z = 0\quad\text{or}\quad A x_1  + B x_2 - \zeta = 0,\\
   \hspace*{5.725em} C x_1 z - \mzs\left(D x_1 + \zeta^{\bbox}\right) = 0;
\end{cases}\label{TheSystem}
\end{gather}
with
\begin{align}
 y_1   = \frac{C}{\mu} x_1 z, \quad 
 y_2   = \frac{\alpha }{\gamma_2} x_2 z, \quad 
  \mys  = \frac{x_1}{\gamma^{\bbox}_1}(C z + \malps\mzs),
  \label{expl-infected-cells}
\end{align}
and
\begin{align}
    A  = \alpha\left(\frac{\gamma_1\kappa_1}{\mu +\gamma_1} - \nu\right), \quad
    B = \alpha(\kappa_2 - \nu), \quad 
    C  = \frac{\alpha  \mu }{\mu +\gamma_1}, \quad 
    D  = \malps(\nu^{\bbox} - 1).
\end{align}
%
%
Observe that since all nondimensionalized parameters are positive, then $C>0$. Also, since $\nu^{\bbox} < 1$ and $\nu < \kappa_{2}$, it follows that $D<0$ and $B>0$. Finally, since $\mu$ is small, then  $\gamma_1\kappa_1/(\mu +\gamma_1) \sim \kappa_1 > \nu$, and thus we may assume that $A >0$.

Solving \eqref{TheSystem} and using \eqref{expl-infected-cells} leads to the nine equilibrium solutions presented in Table \ref{TheCPTable}; full details of the calculations are provided in Appendix \ref{app:CP}.
It is easy to check that the eight conditions, listed in Column~3 of Table~\ref{TheCPTable},
are mutually excluding and contain all possibilities.
The variables $x_1,\ x_2, \ z,$ and $\mzs$ that correspond to the solution of System~\eqref{TheSystem} are listed in the block consisting of Columns~4--7. The block with light yellow background in Columns~8--10 lists the corresponding values of $y_1,\ y_2,$ and $\mys$ obtained from \eqref{expl-infected-cells}. Since all nondimensionalized parameters, together with $A, B, C$, and $D$, are different from zero, it follows that the nine points obtained in Table 1 are well defined and are all the equilibrium points of Eqs.~\eqref{odes_twocell_spec_gen_nd}.
The components of the equilibria of System~\eqref{odes_twocell_spec_gen_nd} need to be non-negative to be biologically meaningful. Ascertaining when this is the case, depending on the values of the nondimensionalized parameters, is not a simple exercise. Thus, this will be done implicitly during the discussion of the stability of equilibria at the numerical level.
\afterpage{
\begin{landscape}
\begin{table}
\begin{center}
\begin{tikzpicture}
\draw [fill=LightCyan, LightCyan] (3.65,0.14) rectangle (8.45,7.75); 
\draw [fill=Palla, Palla] (16.85, 0.14) rectangle (23.18,7.75); 
\node[anchor=south west] at (0,0) {%
  \small\setlength{\tabcolsep}{1.3ex}
  \begin{tabular}{@{}cl
                  r@{\hskip5pt}r@{\hskip5pt}r
                  llll
                  lll
                  @{}
  }\toprule
     \textbf{Critical} & 
     \multirow{2}{*}{\textbf{Name}} &
     \multicolumn{3}{c}{\multirow{2}{*}{\textbf{Conditions}}} & 
     \multirow{2}{*}{{\boldmath{$x_1^*$}}} & 
     \multirow{2}{*}{{\boldmath{$x_2^*$}}} & 
     \multirow{2}{*}{{\boldmath{$z^*$}}} & 
     \multirow{2}{*}{{\boldmath{$\mzs{}^*$}}} &
     \multirow{2}{*}{{\boldmath{$y_1^*$}}} & 
     \multirow{2}{*}{{\boldmath{$y_2^*$}}} & 
     \multirow{2}{*}{{\boldmath{$\mys{}^*$}}} \\
     \textbf{Point} &  & &&& & & & & & & \\ \midrule
     {\boldmath{$v_0$}} & trivial & 
     \multicolumn{3}{r}{$x_1 = 0$ and $x_2 = 0$} & 
     0 & 0 & 0 & 0 & 
     0 & 0 & 0 \\ \midrule
     {\boldmath{$v_1$}} & v-free-1 & 
     \multirow{4}{*}{$x_1 \ne 0;$ $x_2 = 0$} & \multirow{2}{*}{$z = 0$} & $\mzs = 0$ & 
     1 & 0 & 0 & 0 & 
     0 & 0 & 0 \\
     {\boldmath{$v_2$}} & gen-free-1 & 
       & & $\mzs \neq 0$ & 
     $-\frac{\zeta^{\bbox}}{D}$ & 0 & 0 & $\frac{1}{\malps} + \frac{\zeta^{\bbox}}{\malps D}$ & 
     0 & 0 & $-\frac{\zeta^{\bbox}\left(\zeta^{\bbox}+D\right)}{D^2 \gamma^{\bbox}_1}$ \\
     {\boldmath{$v_3$}} & coex-1 & 
       & $z \neq 0$ & &
     $\frac{\zeta}{A}$ & 0 & $\frac{1 - x_1^* - \malps\mzs{}^*}{\alpha}$ 
        & $\frac{1-x^*_1}{\tfrac{\alpha\zeta^{\bbox}}{C x^*_1} + \malps + \tfrac{D \alpha}{C}}$ & 
     $\frac{C}{\mu} x_1^* z^*$ & 0
             & $\frac{x_1^*}{\gamma^{\bbox}_1}(C z^* + \malps\mzs{}^*)$ \\ \midrule
     {\boldmath{$v_4$}} & v-free-2 & 
     \multicolumn{2}{r}{\multirow{2.3}{*}{$x_1 = 0;$ $x_2 \neq 0$}} & $z = 0$ & 
     0 & $\tfrac{\beta_1}{\beta_2}$ & 0 & 0 & 
     0 & 0 & 0 \\
     {\boldmath{$v_5$}} & spec-free & 
       & & $z \neq 0$ & 
     0 & $\frac{\zeta}{B}$ & $\frac{\beta_1}{\alpha} - \frac{\beta_2 \zeta}{\alpha B}$ & 0 & 
     0 & $\frac{\zeta(\beta_1 B - \beta_2 \zeta)}{B^2\gamma_2}$ & 0 \\ \midrule
     {\boldmath{$v_6$}} & gen-free-2 & 
     \multicolumn{2}{r}{\multirow{3.25}{*}{$x_1 \neq 0;$ $x_2 \neq 0$}} & $z = 0$ & 
     $-\frac{\zeta^{\bbox}}{D}$ 
        & $\frac{\beta_1}{\beta_2} + \frac{\zeta^{\bbox}}{D}$ 
        & 0 & $\frac{\beta_2-\beta_1}{\malps\beta_2}$ & 
     0 & 0 & $\frac{(\beta_1-\beta_2)\zeta^{\bbox}}{D\beta_2 \gamma^{\bbox}_1}$\\
     {\boldmath{$v_7$}} & coex-2 & 
       & & \multirow{2}{*}{$z \neq 0$} & 
     $x_1^+$ & $\frac{\zeta - A x_1^+}{B}$
             &  $\frac{\beta_1 - \beta_2\left(x_1^+ + x_2^*\right)}{\alpha}$ 
             &  $\frac{1-\beta_1 + (\beta_2-1)\left(x_1^+ + x_2^*\right)}{\malps}$ & 
     $\frac{C}{\mu} x_1^+ z^*$
             & $\frac{\alpha }{\gamma_2} x_2^* z^*$
             & $\frac{x_1^+}{\gamma^{\bbox}_1}(C z^* + \malps\mzs{}^*)$ \\            
     {\boldmath{$v_8$}} & coex-3 & 
       & & & 
     $x_1^-$ & $\frac{\zeta - A x_1^-}{B}$ 
             & $\frac{\beta_1 - \beta_2\left(x_1^- + x_2^*\right)}{\alpha}$ 
             & $\frac{1-\beta_1 + (\beta_2-1)\left(x_1^- + x_2^*\right)}{\malps}$ &
     $\frac{C}{\mu} x_1^- z^*$
             & $\frac{\alpha }{\gamma_2} x_2^* z^*$
             & $\frac{x_1^-}{\gamma^{\bbox}_1}(C z^* + \malps\mzs{}^*)$ \\
\bottomrule      
\end{tabular}};
\end{tikzpicture}
\end{center}
\captionsetup{singlelinecheck=off}
\caption[.]{\footnotesize
The nine equilibrium points of Eqs.~\eqref{odes_twocell_spec_gen_nd}.  The equilibria in each row is computed by imposing the conditions in Column~3
(with light blue background)
to System~\eqref{TheSystem}, and solving the restricted system with respect to 
the remaining variables. 
The values of $y_1^*,\ y_2^*,$ and $\mys{}^*$ in Columns~8, 9 and 10
(with light yellow background)
are obtained by direct substitution of $x_1^*,\ x_2^*,\ z^*$ and 
$\protect\mzs{}^*$ of the corresponding critical point in \eqref{expl-infected-cells}.
The eight conditions listed in Column~3 are mutually excluding and complete.
$x_1^{\pm}$ denote the two solutions of the quadratic equation 
$\phi_2 x_1^2 + \phi_1 x_1 + \phi_0 = 0,$ where
\begin{align*}
\phi_2 & = (A-B)\Bigl(\bigl(D\alpha + C\malps\bigr)\beta_2 - D\alpha\Bigr), \\
\phi_1 & = (A-B)\alpha\zeta^{\bbox}\bigl(\beta_2	- 1\bigr) + 
           \bigl(D\alpha + C\malps\bigr)\bigl(B\beta_1 - \zeta\beta_2\bigr) + D \alpha(\zeta-B),\text{ and}\\
\phi_0 & = \alpha\zeta^{\bbox}\Bigl(\bigl(1-\beta_2\bigr)\zeta + B\bigl(\beta_1 - 1\bigr)\Bigr).
\end{align*}
\textbf{Notation:} For simplicity we denote respectively by 
$x_1^*,\ x_2^*,\ z^*$ and $\protect\mzs{}^*$ the values of the variables 
$x_1,\ x_2,\ z$ and $\protect\mzs$ 
of a critical point (of course in the same row; that is, for the same conditions).
}\label{TheCPTable}
\end{table}
\end{landscape}
} 

There are two equilibria, $\bm{v}_1$ and $\bm{v}_4$ that describe 
virus-free states and are labelled as ``v-free-1'' and ``v-free-2'',
respectively. 
The two equilibria $\bm{v}_2$ and $\bm{v}_6$ involve the extinction of the 
generalist virus and are labelled as ``gen-free-1'' and ``gen-free-2''. 
The equilibrium point $\bm{v}_5$ corresponds to the extinction
of the specialist and is termed the ``spec-free'' state. 
The final three equilibria, $\bm{v}_3,\ \bm{v}_7,$ and $\bm{v}_8,$
describe various states where the generalist and specialist viral populations 
coexist and they are named ``coex-1'', ``coex-2'' and ``coex-3'' (see Table~\ref{TheCPTable}).

\subsection{Linear stability and codimension-one bifurcations of equilibria}

We now analyze the linear stability of each equilibrium point listed in Table \ref{TheCPTable} and identify codimension-one bifurcations. This is done analytically whenever possible. However, when there is a need to use numerical parameter values, we will set:
\begin{align}
\label{params7d}
\beta_1 = 1.50, \qquad \beta_2 = 2.00, \qquad \mu = 0.10, \qquad \mgammas = \mgammag = \mgammagg = 0.25, \\
\mkg = \mkgg = 1.00, \qquad \mnus = \mnug = 0.50,\qquad \mzets = \mzetg = 0.22, 
\end{align}
and treat the infection rates $\malps$ and $\malpg$ as bifurcation parameters.  
The virulence, burst size, multiplicity of infection, and mortality rate of the viral strains are assumed to be the same for the generalist and the specialist.  For a special case, asymmetric parameters were studied in~\cite{nurtay2019theoretical}, where the asymmetry in the parameters did not introduce new dynamics.

A straightforward calculation shows that the trivial equilibrium $\bm{v}_0$ provides the following vector of eigenvalues:
\begin{equation}
\Lambda_0 = \left(1, -\mzets, -\mzetg, -\mgammas, -\mgammag-\mu, -\mgammagg, \beta_2 \right)^T.
\end{equation}
The first eigenvalue, $1$, is constant due to our choice of nondimensionalization and forces the trivial solution to be a saddle (with a two-dimensional unstable manifold) for all biologically meaningful values of the parameters.

At the first virus-free state $\bm v_1$ (v-free-1), we obtain two negative eigenvalues, $-1$ and $-\mgammagg$.  The third eigenvalue, $\beta_1 - \beta_2$, leads to a case where the growth rates of uninfected cells must be compared. For $\beta_2 > \beta_1$,  $\bm v_1$ can be stable; while for $\beta_2 < \beta_1$,  $\bm v_1$ is unstable.  The next four eigenvalues $\lambda$ are determined from two quadratic equations given by
\begin{subequations}
\begin{align}
Q_1^{(1)}(\lambda) &= \lambda^2+(\malpg\mnug+\mgammag+\mu+\mzetg)\lambda-((\mkg-\mnug)\mgammag-\mu\,\mnug)\malpg+\mzetg(\mgammag+\mu) = 0,\\
Q_2^{(1)}(\lambda) &= \lambda^2+(\malps\mnus+\mgammas+\mzets)\lambda+(\malps(\mnus-1)+\mzets)\mgammas = 0.
\end{align}
\end{subequations}
Rather than analyzing the eigenvalues that arise from solving these quadratic equations, we use these expressions to determine the locations of bifurcation points by forcing the eigenvalue to be zero.  Purely imaginary complex-conjugate eigenvalues are not possible due to the positivity of the linear terms in these expressions.  From $Q_1^{(1)}(0) = 0$ and $Q_2^{(1)}(0) = 0$, we find branches of transcritical bifurcations given, respectively, by
\begin{align}
T_{13} &:= \left\{\malpg = \frac{(\mgammag+\mu)\mzetg}{(\mkg-\mnug)\mgammag - \mu\mnug}\right\}.
\label{eqn:alpha13} \\
T_{12} &:= \left\{\malps = \frac{\mzets}{1-\mnus}\right\}. \label{eqn:alpha12}
\end{align} 
Along the curve given by \eqref{eqn:alpha13}, the equilibria $\bm{v}_1$ and $\bm{v}_3$ collide; along \eqref{eqn:alpha12}, $\bm{v}_1$ and $\bm{v}_2$ collide.  The bifurcation $T_{12}$ represents the emergence of the specialist strain at a critical value of its infectivity $\malps$.  The bifurcation $T_{13}$ is similar but in this case both the specialist and the generalist emerge.  More specifically, the generalist appears when its infection rate $\malpg$ reaches a critical value given by \eqref{eqn:alpha13}.  Due to mutation, the existence of the generalist also gives rise to the specialist.


For the first generalist-free state $\bm v_2$ (gen-free-1), the characteristic polynomial for the eigenvalues $\lambda$ can be written as $L_1^{(2)}(\lambda) L_2^{(2)}(\lambda) Q_1^{(2)}(\lambda) C_1^{(2)}(\lambda) = 0$ where $L_i^{(2)}$, $Q_1^{(2)}$, and $C_1^{(2)}$ denote factors which are linear, quadratic, and cubic in $\lambda$.  In terms of the infection rates, there are three parameter combinations that lead to transcritical bifurcations.  The first coincides with \eqref{eqn:alpha12}.  
The second and third are given by
\begin{align}
T_{26} &:= \left\{\malps = \frac{\beta_2}{\beta_1}\frac{\mzets}{1-\mnus}\right\}, \label{eqn:alpha26} \\
T_{23} &:= 
\left\{\malps = \left[\frac{\mzets}{\mzetg}\frac{(\mkg-\mnug)\mgammag -\mu\mnug}{(\mgammag+\mu)(1-\mnus)}\right]\malpg \right\}, \label{eqn:alpha23}
\end{align}
corresponding to locations where $\bm{v}_2$ collides with $\bm{v}_6$ and $\bm{v}_3$, respectively.  The bifurcation $T_{26}$ describes the situation whereby the infection rate of the specialist becomes so large that there is a substantial decrease in the population of the host $x_1$, allowing the second host, $x_2$, to emerge.  In the case of $T_{23}$, this bifurcation describes the specialist outcompeting the generalist as the infection rate $\malps$ increases.  
The cubic factor $C_1^{(2)}$ gives rise to a Hopf bifurcation curve, denoted by $H_2$, that occurs for a fixed value of $\malps$ that is independent of $\malpg$.  
The quadratic factor $Q_1^{(2)}$ can also produce a Hopf bifurcation. However, it lies outside of the biologically relevant parameter space, \emph{i.e.}\ it occurs for negative infection rates.


The first coexistence state $\bm v_3$ (coex-1) results in a complicated characteristic polynomial. However, it is still possible to recover three analytical parameter combinations which correspond to transcritical bifurcations.  The first two are given by \eqref{eqn:alpha13} and \eqref{eqn:alpha23}. The third combination denotes a curve where $\bm{v}_3$ is equal to $\bm{v}_7$ and is given by
\begin{subequations}
\label{eqn:alpha37}
\begin{align}
T_{37} = \left\{
\malps = \frac{\mzets}{\mzetg}\frac{K^2_0(\beta_1-1)\malpg - K_0\,\mzetg(\beta_2-1)(\mgammag+\mu)}{K_0\,K_1\,\malpg + K_2}\right\},
\end{align}
where
\begin{align}
K_0 &= (\mkg-\mnug)\mgammag - \mu\mnug,\\
K_1 &= (\beta_1-1)(1-\mnus)\mgammag+(1+(\beta_1-1)\,\mnus)\mu,\\
K_2 &= (\mgammag+\mu)((1+\mnus(\beta_2-1))\mu-(1-\mnus)(\beta_2-1)\mgammag)\mzetg.
\end{align}  
\end{subequations}
The bifurcation $T_{37}$ describes a similar scenario as $T_{26}$, whereby the high infectivity rate of the specialist reduces the population of $x_1$, allowing $x_2$ to emerge.  However, the generalist is now present when this happens, leading to the appearance of the infected cell population $\mygg$ after the bifurcation occurs.

The second virus-free state $\bm{v}_4$ (v-free-2) has four negative eigenvalues:
$-\mzets$, $-\mgammag-\mu$, $-\mgammas$, and $-\beta_1$.  The fifth eigenvalue, $1 - \beta_1 / \beta_2$, is only negative if $\beta_1 > \beta_2$. Thus, we may have an exchange of stability between the v-free-1 and v-free-2 states as the relative values of $\beta_1$ and $\beta_2$ change.  Recall that v-free-1 could become stable only if $\beta_2 > \beta_1$.  This change in stability reflects the competitive exclusion principle: as both $x_1$ and $x_2$ are competing for the same resource, only the fittest cell type will survive.  The last two eigenvalues can be determined from 
\begin{align}
Q^{(4)}(\lambda) = \beta_2\lambda^2+[(\mzetg+\mgammagg)\beta_2+\mnug\,\beta_1\malpg]\lambda+\mgammagg\,[\mzetg\,\beta_2-\beta_1(\mkgg-\mnug)\malpg].
\label{eqn:v4_Q}
\end{align}
The form of \eqref{eqn:v4_Q} shows that Hopf bifurcations cannot occur.  However, there is a curve of transcritical bifurcations given by
\begin{equation}
T_{45} = \left\{\malpg = \frac{\beta_2}{\beta_1}\frac{\mzetg}{\mkgg-\mnug}\right\},
\label{eqn:alpha45}
\end{equation}
where $\bm{v}_4$ and $\bm{v}_5$ intersect.  The bifurcation $T_{45}$ is analogous to $T_{12}$ and describes the emergence of the generalist strain when its infectivity exceeds a critical value. 

The specialist-free state $\bm v_5$ (spec-free) depicts the only case of the entire viral population consisting solely of the generalist strain.  There are three eigenvalues $-\mzets$, $ -\mgammag - \mu$, $-\mgammas$ which are always negative.  The remaining four are determined by the solutions of cubic and linear equations:
\begin{subequations}
\begin{align}
&\malpg\,(\mkgg-\mnug)^2\lambda^3+(\mkgg-\mnug)((\malpg\,\mkgg+\beta_2)\mzetg+\mgammagg\malpg\,(\mkgg-\mnug))\lambda^2 \nonumber\\
&\quad +(\beta_2(\mkgg+\mnug)\mzetg-(\mkgg-\mnug)(\malpg\,\beta_1\mnug-\beta_2\mgammagg))\mzetg\lambda+(\mkgg-\mnug)(\beta_1\malpg\,(\mkgg-\mnug)-\mzetg\beta_2)\mgammagg\mzetg = 0, \label{eqn:v5_cubic}\\
&\malpg\,(\mkgg-\mnug)\,\lambda+(\mkgg-\mnug)(\beta_1-1)\malpg-\mzetg(\beta_2-1) = 0. \label{eqn:v5_linear}
\end{align}
\end{subequations}
Setting $\lambda = 0$ in \eqref{eqn:v5_cubic} leads to the transcritical bifurcation $T_{45}$ given by \eqref{eqn:alpha45}.  Setting $\lambda = \pm i \omega$ shows there is a curve of Hopf bifurcations denoted by $H_5$ that occurs along a fixed value of $\malpg$ that is independent of $\malps$.  Finally, substituting $\lambda = 0$ into \eqref{eqn:v5_linear} produces another curve of transcritical bifurcations involving $\bm{v}_5$ and $\bm{v}_7$ given by
\begin{equation}
T_{57} := \left\{\malpg = \frac{\beta_2-1}{\beta_1-1}\frac{\mzetg}{\mkgg-\mnug}\right\}.
\label{eqn:alpha57}
\end{equation}
The bifurcation $T_{57}$ is similar to $T_{23}$; in this case, increasing $\alpha$ beyond that given by \eqref{eqn:alpha57} leads to the generalist outcompeting the specialist.

The generalist-free state $\bm v_6$ (gen-free-2) also leads to a complicated characteristic polynomial that can be factored into the product of a cubic and quartic polynomial in $\lambda$. Nevertheless, it is possible to recover two parameter combinations for zero-eigenvalue bifurcations.  The first coincides with $T_{26}$ in \eqref{eqn:alpha26}.  The second corresponds to a transcritical bifurcation with the coexistence state $\bm{v}_7$ and is given by
\begin{equation}
T_{67} = \left\{\malps = \malpg \, \frac{\beta_2\mzets\,(\mkg\,\mgammag-\mkgg(\mu+\mgammag))}{(1-\mnus)(\mu+\mgammag)(\beta_1\mzetg - \beta_1\malpg(\mkgg-\mnug))}\right\}.
\label{eqn:alpha67}
\end{equation}
$T_{67}$ indicates the threshold where $\bm{v}_6$ loses stability to $\bm{v}_7$ as $\malps$ is increased.  The cubic polynomial can give rise to a curve of Hopf bifurcations but this leads to some components of $\bm{v}_6$ being negative for the parameter values we consider. From the quartic polynomial we can deduce the existence of a curve of Hopf bifurcations, denoted as $H_6$, that occur for fixed values of $\malps$ independent of $\malpg$.  

The coexistence equilibria $\bm{v}_7$ and $\bm{v}_8$ are determined by solving a complicated quadratic polynomial.  Consequently, it is not practical to study their linear stability through analytical means.  However, using numerical continuation, we find that this pair of equilibria emerge \emph{via} a curve of saddle-node bifurcations denoted by LP$_{7}$.  From a biological perspective, this curve represents the theoretical threshold that must be exceeded in order for all seven populations to coexist.  Both $\bm{v}_7$ and $\bm{v}_8$ can undergo Hopf bifurcations along the curves $H_7$ and $H_8$.  However, for the parameters that we consider, only the curve $H_7$ lies in the biologically meaningful region of phase space.



All bifurcation curves can be simultaneously plotted in the $(\malpg,\ \malps$) plane, giving rise to the complex two-dimensional bifurcation diagram shown in Fig.~\ref{fig:bif_diag}. 
In addition to the transcritical, saddle-node, Hopf bifurcations that were identified through the linear stability analysis, there are also curves of saddle-node and transcritical bifurcations of cycles that emanate from codimension-two points, which we now describe. 


\afterpage{
\begin{landscape}
\begin{figure} 
\centering
\captionsetup{width=\linewidth}
\includegraphics[width=1\linewidth]{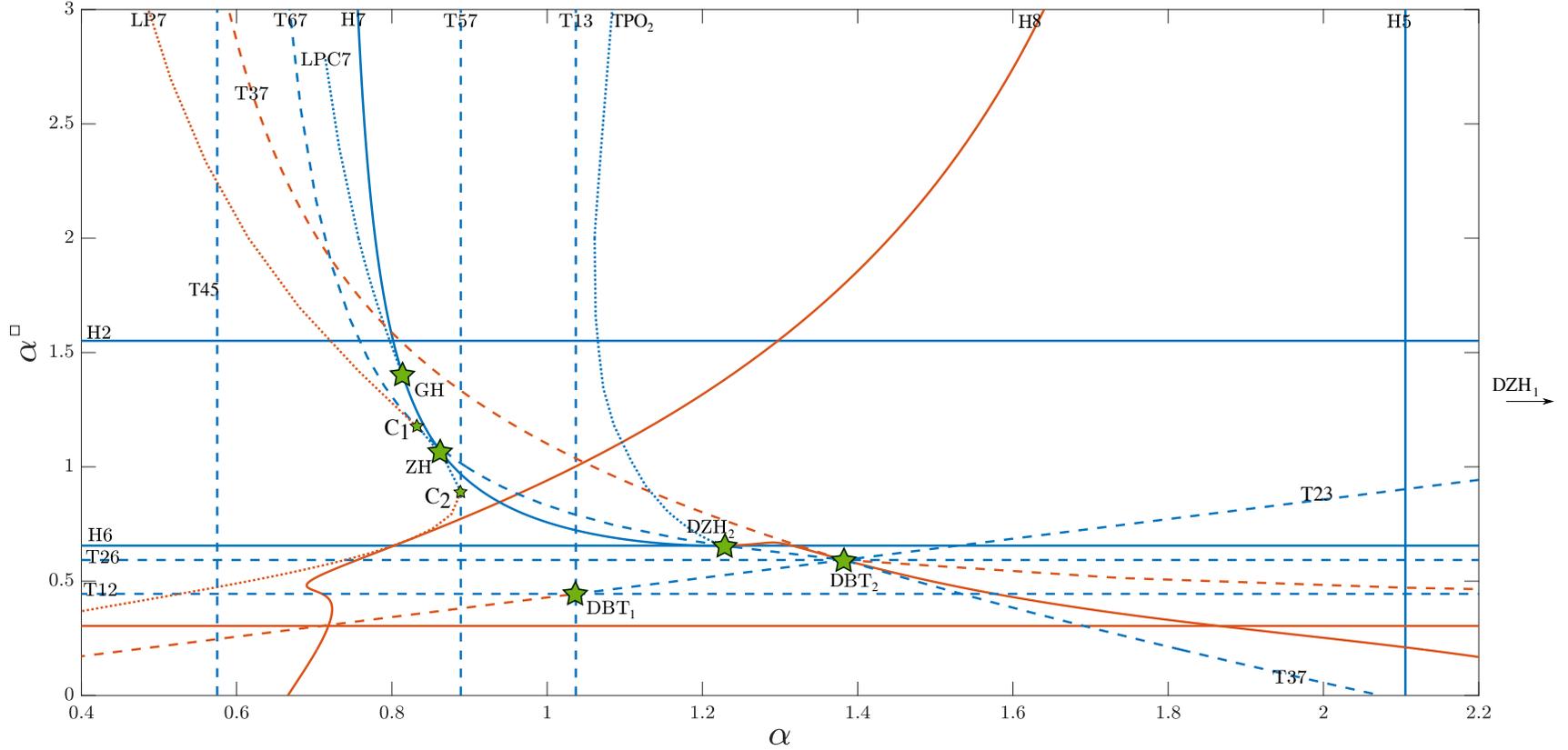}
\caption{\small Two-dimensional bifurcation diagram for system \eqref{odes_twocell_spec_gen_nd} with bifurcation parameters \protect \alpg and \protect \alps. The remaining parameters are fixed as in \eqref{params7d}. Dashed and solid lines correspond to transcritical ($T_{ij}$) and Hopf ($H_j$) bifurcations, respectively. Dotted curves correspond to limit points of equilibria (LP), limit point of cycles (LPC), or transcritical bifurcation of periodic orbits (TPO). Blue denotes biologically relevant solutions (strictly positive values for all components). Orange describes biologically irrelevant solutions with negative components. The green stars denote codimension-two bifurcations, \emph{i.e.}\ generalised Hopf (GH), zero-Hopf (ZH), degenerate zero-Hopf (DZH), degenerate Bogdanov-Takens (DBT) bifurcations, and cusp-like bifurcations~(C).}
\label{fig:bif_diag}
\end{figure}
\end{landscape}
}

\subsection{Identification and classification of codimension-two bifurcations}

The intersection of the transcritical curves $T_{12}$ and $T_{13}$ results in a double transcritical bifurcation.  The Jacobian of the system \eqref{odes_twocell_spec_gen_nd} at this point has a double-zero eigenvalue and therefore it is a Bogdanov-Takens (BT) bifurcation.  By calculating the normal form coefficients following the procedure and notation of Kuznetsov, we find that the bifurcation is degenerate because the normal form coefficient $a_2 = 0$ (see Ref.~\cite{kuznetsov2005practical}).  We thus refer to this as a degenerate BT (DBT) bifurcation and label it as DBT$_1$.  A similar bifurcation was detected in \cite{nurtay2019theoretical}, who studied its properties by seeking series expansions of the eigenvalues near this point.  It was shown that only curves of transcritical bifurcations emanate from this point. The same conclusion is reached in this case as well.

The four transcritical curves $T_{26}$, $T_{23}$, $T_{37}$, and $T_{67}$ given by \eqref{eqn:alpha26}--\eqref{eqn:alpha37} and \eqref{eqn:alpha67} also intersect at a single point, leading to a quadruple transcritical bifurcation, whereby four equilibria, $\bm{v}_2$, $\bm{v}_3$, $\bm{v}_6$, and $\bm{v}_7$ coincide.  At this point, the Jacobian of \eqref{odes_twocell_spec_gen_nd} has a double-zero eigenvalue, indicating a type of BT bifurcation.  However, the Jacobian has two linearly independent eigenvectors associated with the zero eigenvalues, a feature which is not generic to BT bifurcations.  Consequently, the Jordan block associated with these eigenvalues does not match the form typically assumed in developing the normal form theory for standard and degenerate BT points~\cite{kuznetsov2005practical, kuznetsov2013elements} and thus the results from these works do not apply.  Nevertheless, we label this point as DBT$_2$.

There are three zero-Hopf (ZH) bifurcations that occur in the system, two of which are degenerate.  The first degenerate zero-Hopf point lies at the intersection of the Hopf curve $H_2$ and the transcritical curve $T_{23}$ and is labelled DZH$_1$.  The chosen parameter values lead to DZH$_1$ laying outside of the biologically meaningful phase space and therefore we do not consider it further.  The second, DZH$_2$, lies at the intersection of $H_6$ and $T_{67}$.  Following the notation of Kuznetsov~\cite[Sec.~8.5]{kuznetsov2013elements}, the degeneracy arises because the normal form coefficient $G_{011}(0) = 0$.  A similar bifurcation, also found in~\cite{nurtay2019theoretical}, was shown to produce a curve of transcritical bifurcations of periodic orbits (TPOs).  Numerical continuation confirms that the same occurs here.  The curves of TPOs that emanate from DZH$_2$ are shown in Fig.~\ref{fig:bif_diag} and are labelled as TPO$_2$. The non-degenerate ZH bifurcation occurs at the tangential intersection of $H_7$ and a saddle-node bifurcation (limit-point bifurcation) of $\bm{v}_7$ and $\bm{v}_8$, LP$_7$.  Again, following Kuznetsov~\cite[Sec.~8.5]{kuznetsov2013elements}, the normal form coefficients are found to be $s = 1$ and $\theta(0) > 0$, indicating the existence of subcritical Hopf bifurcations, as will be confirmed below.  For this type of ZH point, there is no branch of torus bifurcations. 

A generalized Hopf bifurcation denoted as GH occurs along the Hopf curve $H_7$.  A curve of saddle-node bifurcations of cycles, which we label as LPC$_7$ emerges from the GH point. We will show that the LPC$_7$ curve plays a very important role because it marks the region of parameter space where chaotic dynamics occur. 

Finally, there are codimension-two bifurcations that occur at the intersection of LP$_7$ with $T_{67}$ and $T_{57}$. These are labelled as $C_1$ and $C_2$ on Fig.~\ref{fig:bif_diag}, respectively. Only the latter is biologically relevant because it gives rise to a bistable region of parameter space where either all seven populations coexist and tend to a stationary equilibrium or the specialist dies out and only $x_2$, $\mygg$, and $\mzg$ survive with oscillatory population sizes. 



\section{One-dimensional bifurcation diagrams}
\label{sec:onedim}

Figure~\ref{fig:bif_diag} reveals a complex bifurcation scenario.  To better understand the dynamics of the system in various regions of parameter space,
we construct one-dimensional bifurcation diagrams by allowing one of the infection rates to vary while fixing the other.  This approach is equivalent to taking horizontal and vertical cross-sections in Fig.~\ref{fig:bif_diag}.  To visualize the different equilibria, we plot how their Euclidean norm
\begin{equation}
V_i = ||\bm v_i|| = \sqrt{x^2_1+x^2_2+(\mys)^2+\myg^2+\mygg^2+(\mzs)^2+\mzg^2}
\end{equation}
evolves under parameter variation.  
The definition of the norm might lead to seemingly intersecting curves on a plot, which may happen when the vectors lay on a sphere of the same radius in a seven-dimensional space. 
When plotting periodic orbits we adapt the same norm for maxima and minima of the periodic orbit.  Furthermore, due to the complexity of the bifurcation diagrams, we only show biologically meaningful solutions, \emph{i.e.} those with components that are non-negative. 

\begin{figure}
	\centering
	\subfigure[$\protect \malps = 0.20$]{\includegraphics[width=0.495\linewidth]{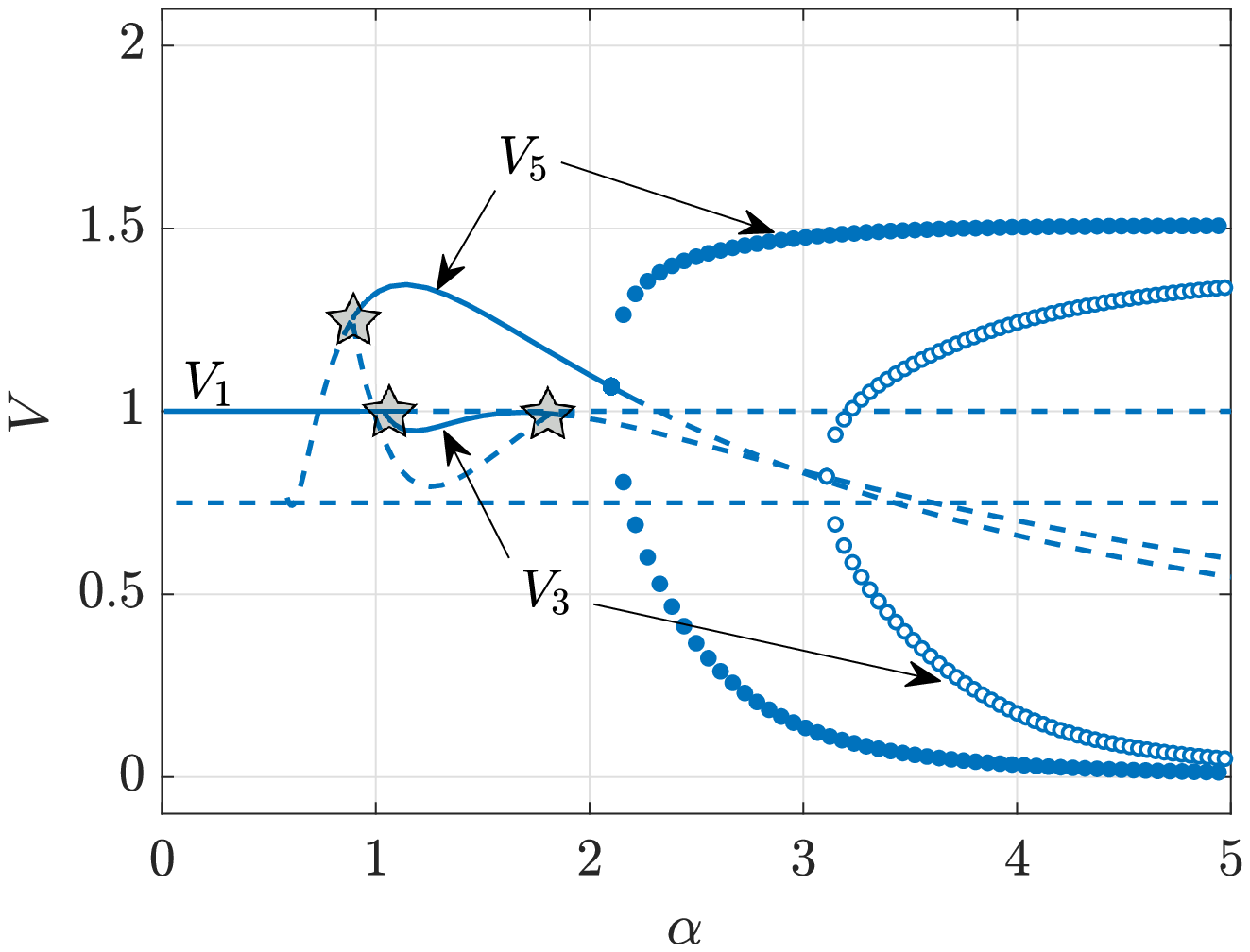}}
	\subfigure[$\protect \malps = 0.53$]{\includegraphics[width=0.495\linewidth]{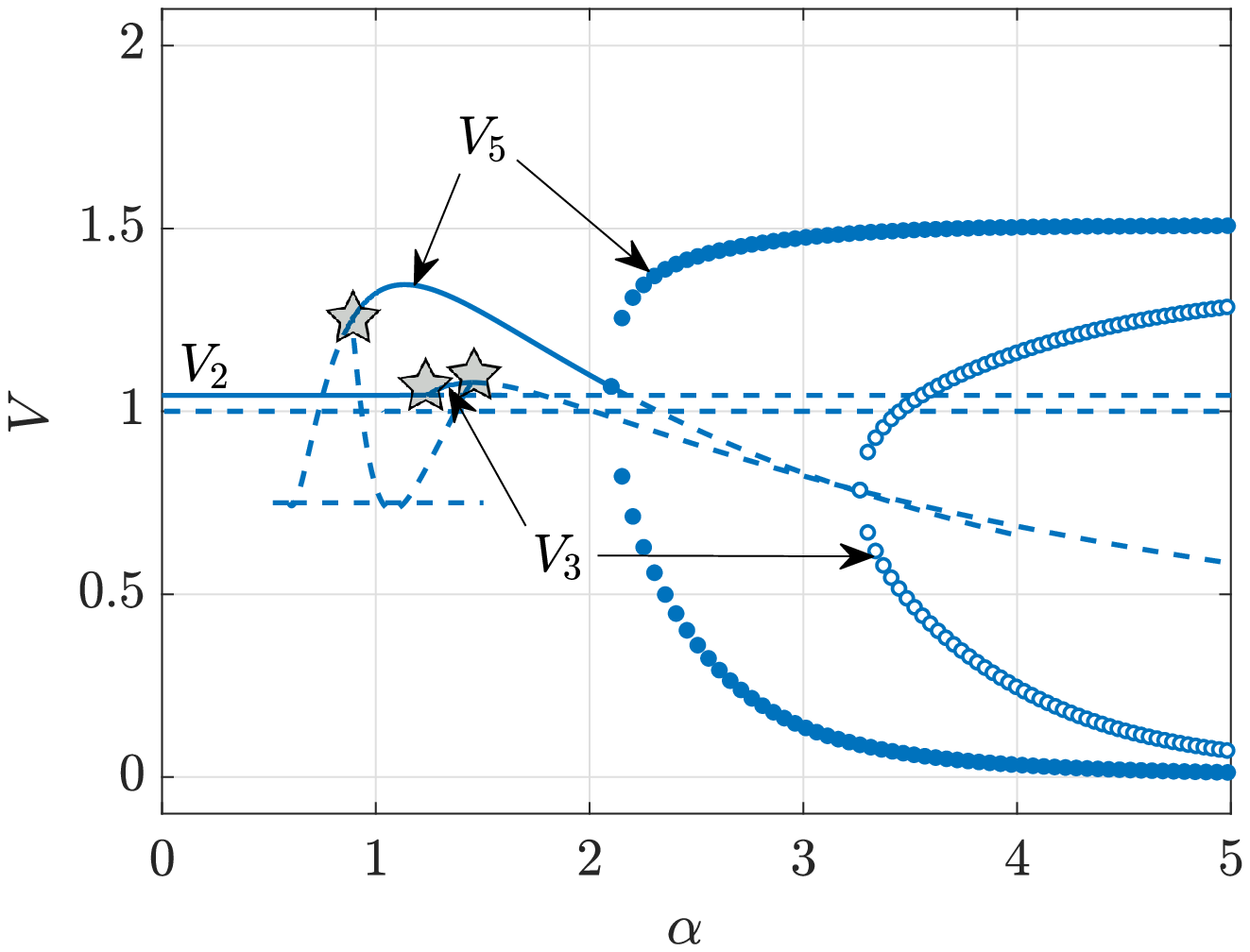}} \\
	\subfigure[$\protect \malps = 0.62$]{\includegraphics[width=0.495\linewidth]{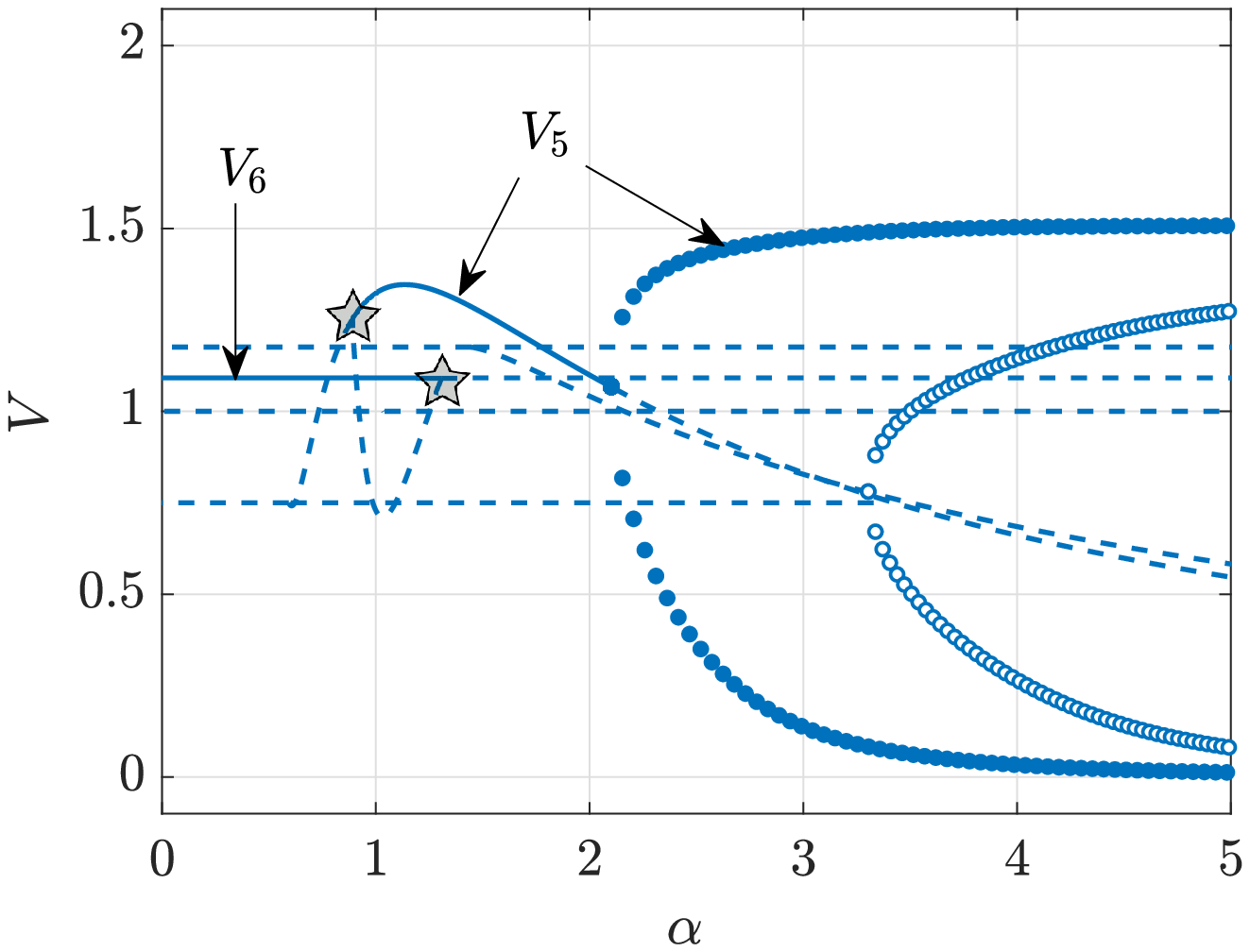}}
	\caption{One-dimensional bifurcation diagrams depicting the norm of all biologically meaningful equilibria, $V_i$, as functions of the generalist infection rate $\protect \malpg$ for fixed values of the specialist infection rate $\protect \malps$ (the values of all other parameters are listed in~\eqref{params7d}).  Solid and dashed lines correspond to stable and unstable equilibria, respectively. Filled and empty markers represent stable and unstable periodic orbits. The grey stars denote the locations of transcritical bifurcations.}
	\label{fig:one_d_1}
\end{figure}

Figure~\ref{fig:one_d_1} (a) illustrates the behavior in the regions below $T_{12}$; it corresponds to a fixed value of $\malps = 0.2$.  For small values of \alpg, the v-free-1 state $\bm v_1$ remains stable until it crosses and exchanges stability with the coex-1 state $\bm v_3$ at $\malpg = 1.04$, corresponding to $T_{13}$.  The stable equilibrium $\bm v_3$ then intersects and exchanges stability with the coex-2 state $\bm v_7$ at $\malpg = 1.82$ (\emph{i.e.}\ $T_{37}$).  The coex-2 state gains stability but moves outside of the biologically meaningful phase space. Following the same unstable (dashed) curve of $\bm v_7$ backwards, we notice that it intersects and changes stability with the spec-free state $\bm v_5$ at $\malpg = 0.89$ (\emph{i.e.}\ $T_{57}$).  Thus, the transcritical bifurcations $T_{57}$ and $T_{37}$ mark the boundaries of a bistable region of parameter space and also lead to hysteresis.  Between the minimum and maximum values of the bistable interval we witness two different types of bistability: the v-free-1 state $\bm v_1$ and spec-free state $\bm v_5$ are simultaneously stable, then the coex-1 state $\bm v_3$ and thee spec-free state $\bm v_5$ are both stable. A supercritical Hopf bifurcation ($H_5$) leads to the creation of stable periodic orbits around $\bm{v}_5$.  The Hopf bifurcation of the coex-1 state $\bm v_3$ ($H_3)$ leads to unstable periodic orbits.  

To explore the regions between the bifurcation curves $T_{12}$ ($\malps = 0.45$) and $T_{26}$ ($\malps = 0.59$), we fix $\malps = 0.53$ to produce the bifurcation diagram shown in Fig.~\ref{fig:one_d_1} (b).
This one-dimensional bifurcation diagram is similar to the one in Fig.~\ref{fig:one_d_1} (a). However, the gen-free-1 state $\bm v_2$ is now stable for small values of \alpg, whereas the v-free-1 state $\bm v_1$ is unstable.  With increasing \alpg, the spec-free state $\bm v_5$ also becomes stable at $\malpg = 0.89$ (\emph{i.e.}\ $T_{57}$).  The bistability between $\bm v_2$ and $\bm v_5$ is interesting because it implies that, depending on the initial conditions, the system may either lose all generalists or all specialists.  That is, although both populations can exist, they cannot coexist. As $\malpg$ is increased to $1.23$, the bistability changes and occurs between the coex-1 state $\bm v_3$ and the spec-free state $\bm v_5$, as before.  However, the region of bistability between these two points is now much smaller.

By increasing $\malps$ to 0.62, the transcritical bifurcation $T_{26}$ is induced, resulting in the gen-free-1 state $\bm{v}_2$ losing stability and the gen-free-2 state $\bm{v}_6$ gaining stability.  The difference between these two states is that the  uninfected cell population of gen-free-1 only consists of $x_1$ whereas the gen-free-2 state allows the coexistence of both $x_1$ and $x_2$.
This increase in $\malps$ also induces the bifurcation DBT$_2$, which eliminates the bistability that exists between the coex-1 state $\bm{v}_3$ and the spec-free state $\bm{v}_5$.  Thus, bistability only occurs between the gen-free-2 state and the spec-free state, which again leads to a scenario where only one of the two viral strains can persist in the long term.

\begin{figure}
	\centering
	\subfigure[]{\includegraphics[width=0.495\linewidth]{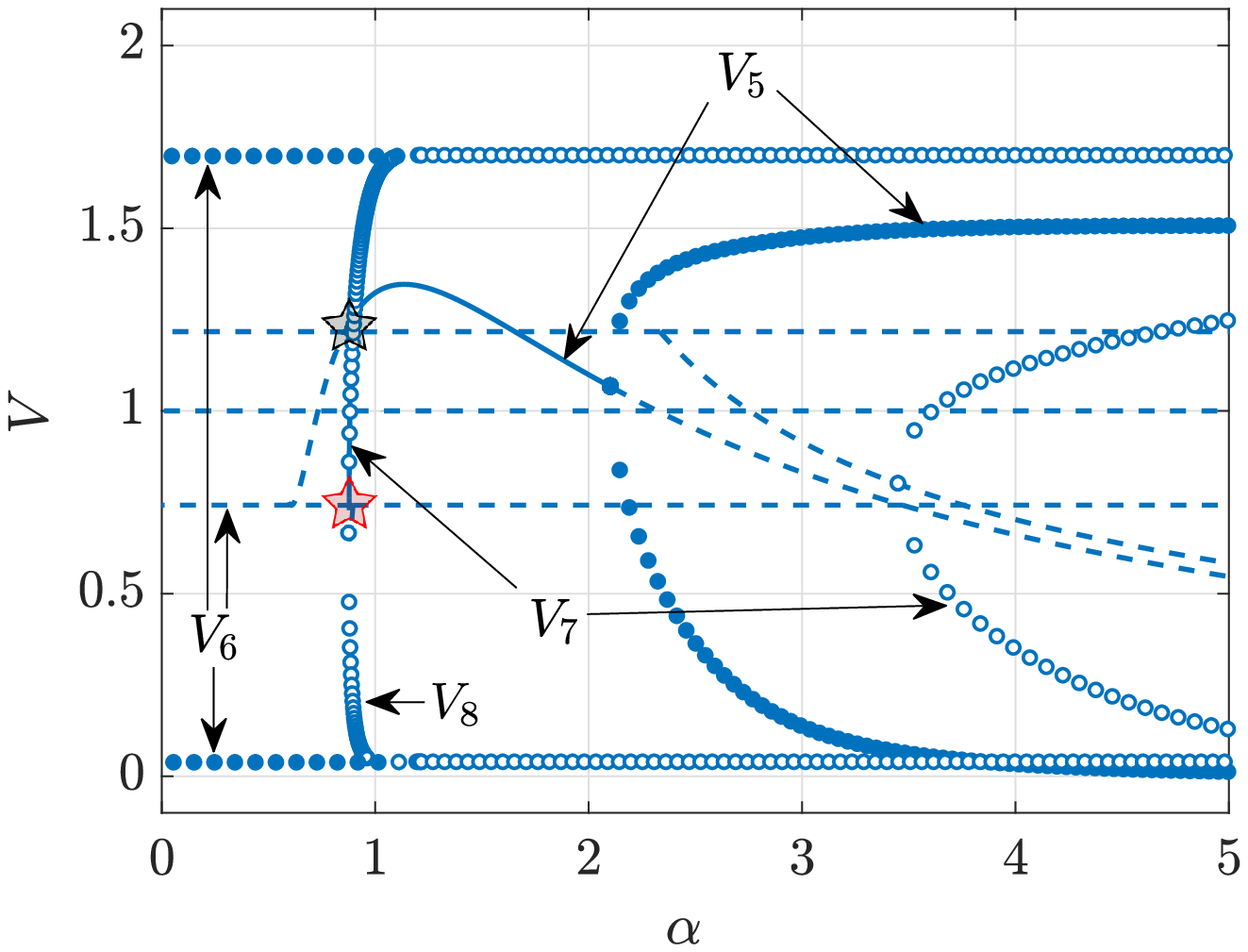}}
	\subfigure[]{\includegraphics[width=0.495\linewidth]{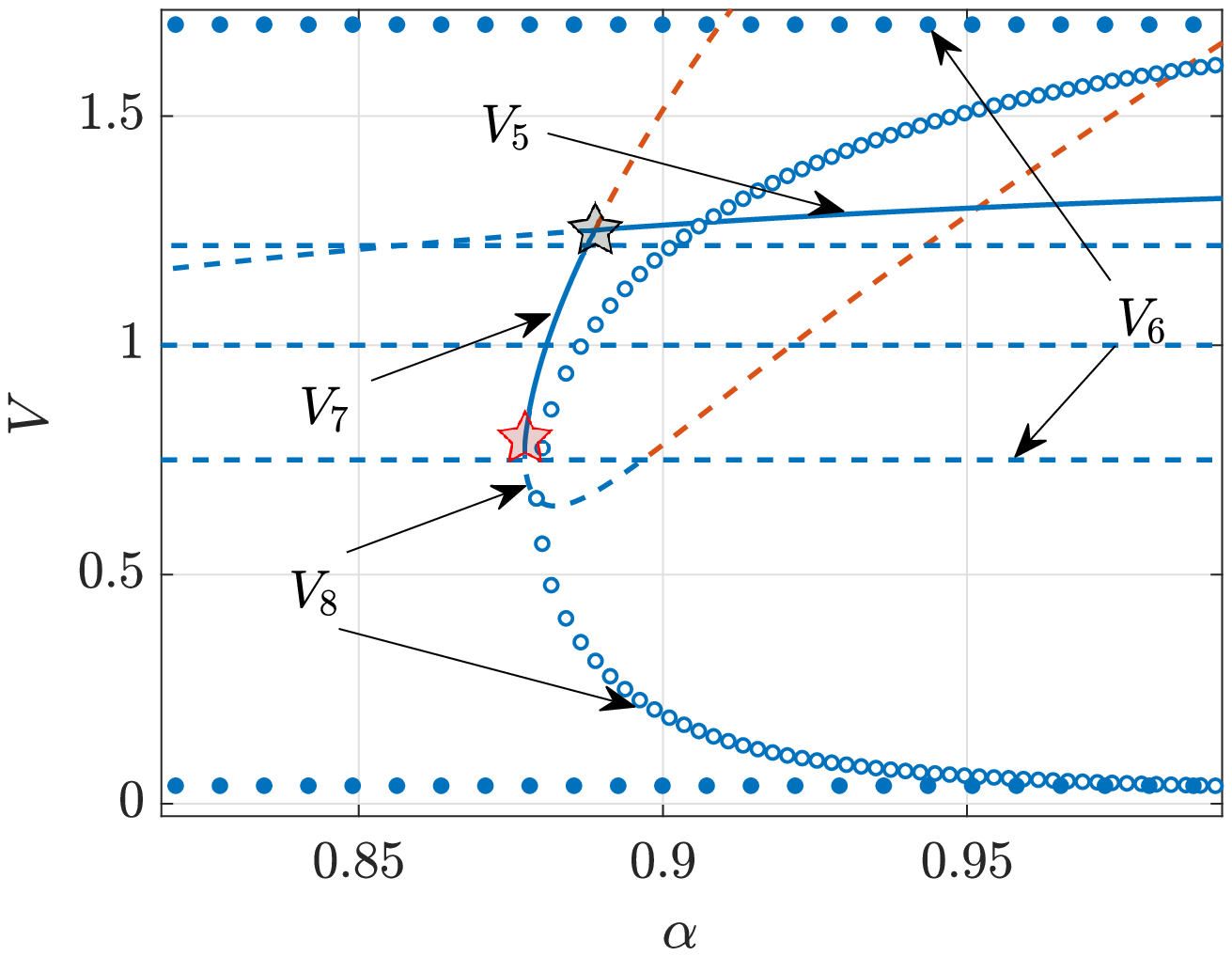}}
	\caption{One-dimensional bifurcation diagrams for a fixed value of the specialist infection rate $\protect \malps = 1.0$.   Panel (b) focuses on the region near $\protect \malpg = 0.9$.  Solid and dashed curves stand for stable and unstable equilibria, respectively. Filled markers and empty markers stand for stable and unstable periodic orbits, respectively. The red star at $\protect \malpg = 0.877$ corresponds to a saddle-node bifurcation LP$_{7}$. The grey star at $\protect \malpg = 0.89$ is the transcritical bifurcation corresponding to $T_{57}$.}
	\label{fig:a_1}
\end{figure}

The one-dimensional bifurcation diagram for $\malps = 1$ is shown in Fig.~\ref{fig:a_1} (a), which illustrates a rich set of dynamics which are seen more clearly when enlarged in Fig.~\ref{fig:a_1} (b).  For sufficiently small values of \alpg, none of the equilibria are stable; instead, there is a stable periodic orbit (PO) which emerges from $\bm{v}_6$ due to the onset of the supercritical Hopf bifurcation (in Fig.~\ref{fig:bif_diag} indicated by $H_6$). 
By increasing $\malpg$ from zero, we observe the saddle-node bifurcation LP$_{7}$, marked with a red star, which creates the coex-2 and coex-3 states at $\malpg = 0.877$.  Here, two branches of equilibria, $\bm v_7$ (stable) and $\bm v_8$ (unstable), appear. The stable branch, $\bm v_7$, undergoes the transcritical bifurcation $T_{57}$ with the unstable spec-free state $\bm v_5$ at $\malpg = 0.89$ which leads to an exchange of stability. The unstable branch, $\bm v_8$, undergoes a subcritical Hopf bifurcation ($H_7$) at $\malpg = 0.879$ which creates unstable POs.  These unstable POs grow and collide with the stable POs that emerged from $\bm{v}_6$, resulting in a TPO bifurcation (in Fig.~\ref{fig:bif_diag} indicated by TPO$_2$). The POs exchange stability at this point but the POs that emerge from $\bm{v}_8$ leave the biologically meaningful phase space.

From the bifurcation diagrams shown in Fig.~\ref{fig:a_1}, we see there are two bistable regions: one between $\malpg = 0.877$ (\emph{i.e.}\ LP$_7$) and $\malpg = 0.89$ (\emph{i.e.}\ $T_{57}$), and the second between $\malpg = 0.89$ (\emph{i.e.}\ $T_{57}$) and $\malpg = 1.11$ (\emph{i.e.}\ TPO$_2$).  
If the infection rate of the generalist lays between LP$_7$ and $T_{57}$, either both viral strains coexist and tend to steady population sizes or the specialist outcompetes the generalist but has a fluctuating (oscillatory) population size.  However, for greater generalist infection rates which lie between $T_{57}$ and TPO$_2$, the viruses cannot coexist in a stationary manner.  Either the generalist outcompetes the specialist and tends to a steady population size or the specialist outcompetes the generalist with a fluctuating population size.  



\begin{figure}
	\centering
	\subfigure[]{\includegraphics[width=0.495\linewidth]{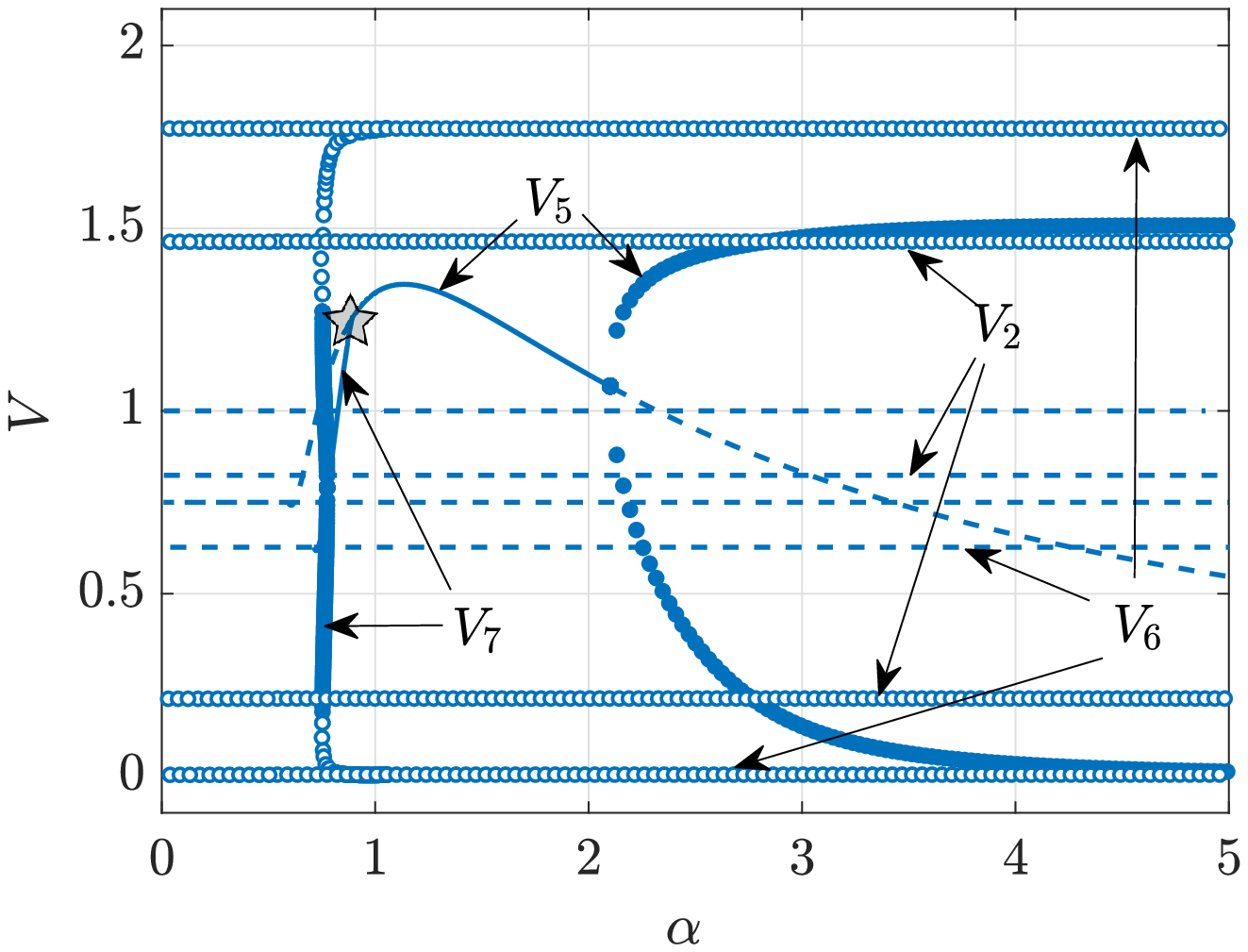}}
    \subfigure[]{\includegraphics[width=0.495\linewidth]{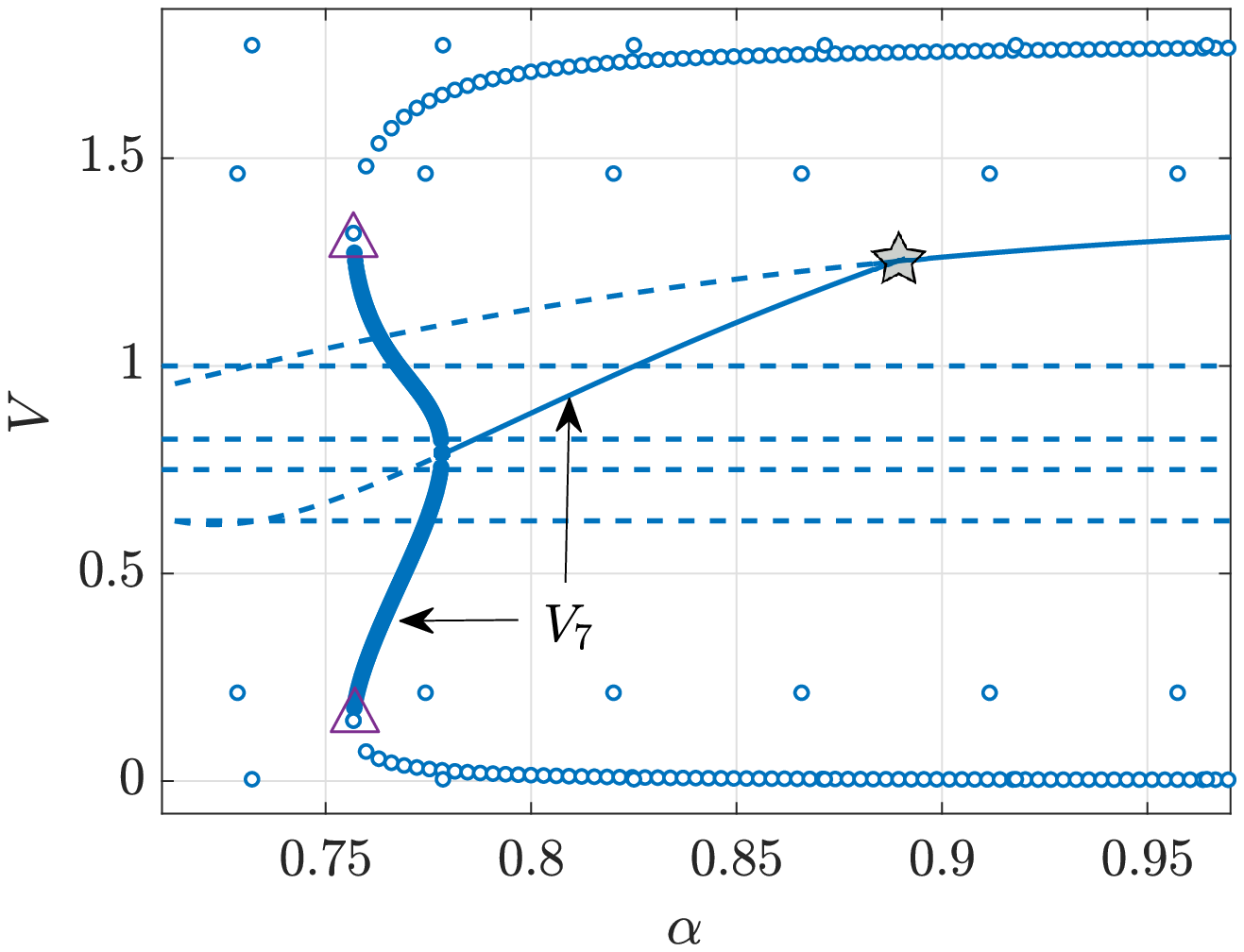}}
	\caption{One-dimensional bifurcation diagram for a fixed value of $\protect \malps = 2$.  Panel (b) focuses on the region near $\protect \malpg = 0.85$.  Solid curves and filled markers stand for stable equilibria and stable periodic orbits, respectively. Dashed curves and empty markers illustrate unstable equilibria and unstable periodic orbits. The grey star at $\protect \malpg = 0.89$ is the transcritical bifurcation $T_{57}$.  The purple triangle marks the limit point of cycles LPC$_1$.}
	\label{fig:a_2}
\end{figure}

An elaborate set of dynamics is observed for greater values of the specialist infection rate. This case is important from a biological perspective because specialization is known to be accompanied by the rise of fitness of a population~\cite{schwartz2012specialisation, leiby2014adaptation}. The one-dimensional bifurcation diagram presented in Fig.~\ref{fig:a_2} corresponds to a value of $\malps = 2$.  For small values of \alpg, there are no stable equilibria nor stable POs. As will be discussed in Section \ref{sec:chaos}, a strange chaotic attractor exists in this region.  
Due to the onset of the ZH bifurcation followed by the GH bifurcation, the subcritical Hopf bifurcation $H_7$ shown in Fig.~\ref{fig:a_1} moves from the branch of $\bm{v}_8$ to the branch of $\bm{v}_7$ and becomes supercritical.
The stable PO that emerge from this Hopf bifurcation grow in amplitude with decreasing $\malpg$ until they destructively collide with an unstable PO at an LPC bifurcation.  This growth corresponds to increases in the populations of the healthy cells ($x_1$ and $x_2$) as well as the specialist ($\mys$, $\myg$, and $\mzs$), which fluctuate with greater amplitudes due to the diminishing population of the generalist. 
The LPC bifurcation marks the point at which stable solutions no longer exist and the onset of chaotic dynamics.  Due to the loss of stability of the POs that emerge from $\bm{v}_6$, there are no bistable regions for this value (or greater values) of the specialist infection rate.


\section{Probabilistic stability diagram} 
\label{sec:prob}
\begin{figure}
	\centering
	\includegraphics[width=1\linewidth]{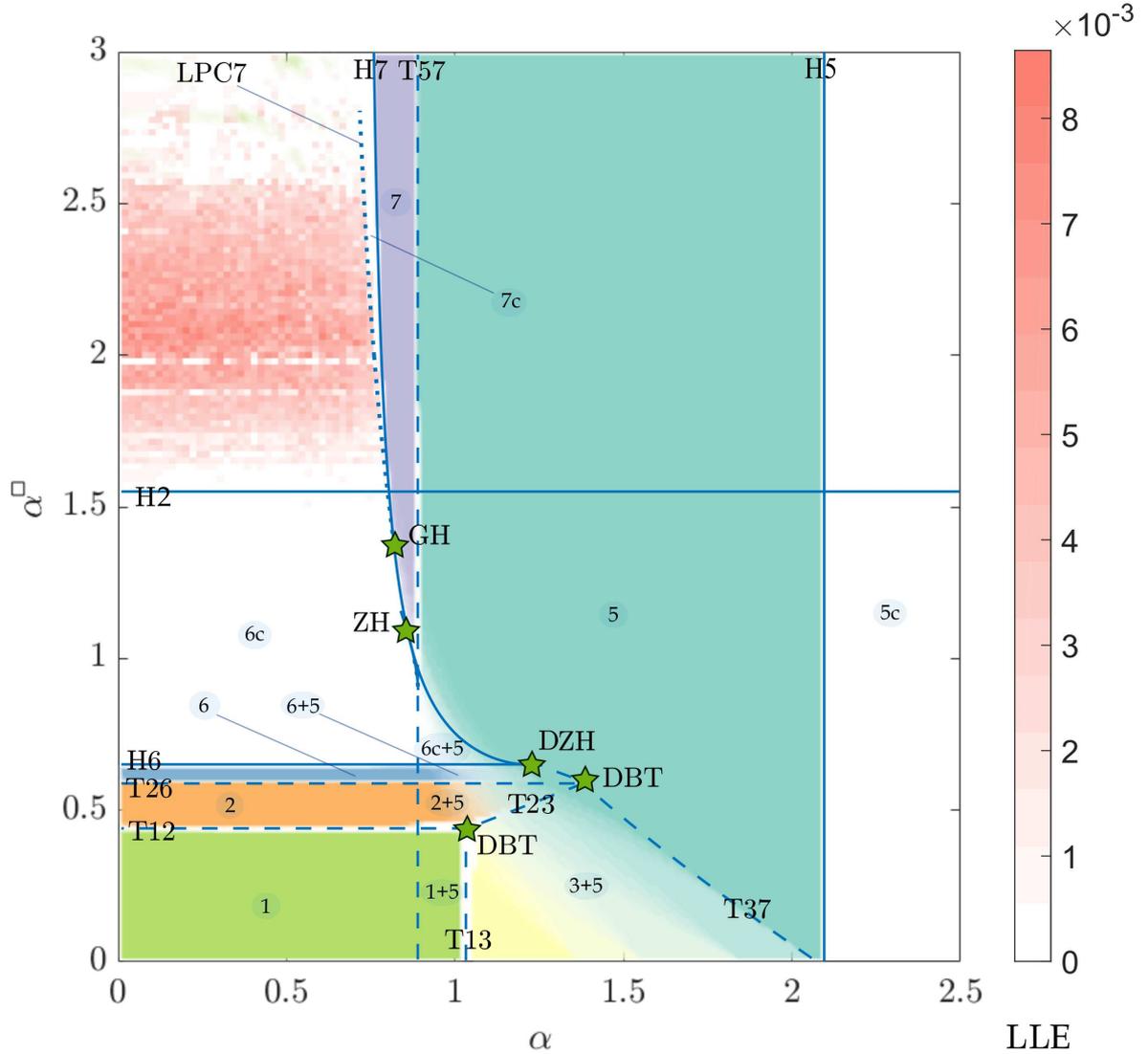}
	\caption{Stability diagram for equilibria $\bm v_i$, $i = 1, \ldots,  7$; regions are marked with indices of the stable equilibria. Intensity of all colors (except red) in regions illustrate the probability of the system stabilizing at the given equilibria when started from different initial conditions (only $\protect \mzs(0)$  and $\protect \mzg(0)$ were sampled). Regions of bistability are indicated with addition, \emph{e.g.} 3+5 when equilibria $\bm{v}_3$ and $\bm{v}_5$ are acheived from different initial conditions. White areas correspond to stability of periodic orbits. Positive largest Lyapunov exponents (LLE, in red) illustrate the region corresponding to chaotic behavior. Dashed, solid, and dotted curves correspond to transcritical bifurcation, Hopf bifurcation, and saddle-node bifurcation of periodic orbits, respectively.}
	\label{fig:allbif}
\end{figure}

The one-dimensional bifurcation diagrams revealed several regions of bistability.  We use a probabilistic approach to explore the basins of attraction and to ascertain the likelihood of the system converging to the various stable states that exist for a given set of parameters.  More specifically, we fix the initial conditions for the susceptible cell populations at the maximum capacity of the environment in the absence of virus, \emph{i.e.}\ $x_1(0)=1.0$ and $x_2(0)=0.75$ for the parameters in~\eqref{params7d}.  The initial viral loads $\mzg(0)$ and $\mzs(0)$ are systematically taken from every grid point of a unit square that has been uniformly discretized into $N\times N$ points (we take $N=5$ to ensure computational feasibility). Consequently, for every fixed set of parameters, $N^2$ simulations are conducted with different initial conditions of viral loads. The probability of becoming attracted to equilibrium $\bm v_i$ is calculated as the number of times the numerical solution converges to that state per $N^2$ runs.  The probabilities of arriving at each possible stable state are mapped in Fig.~\ref{fig:allbif} with different colors.  Higher probabilities are described by increasing the intensity of the color.  White denotes regions where there are stable periodic solutions.  Red denotes chaotic regions of parameter space with the intensity being linked to the value of the largest Lyapunov exponent, as described in the next section. 
Regions marked by $(i)$ or $(i+j)$ are those in which $\bm{v}_i$ is stable (monostable regions) or $\bm{v}_i$ and $\bm{v}_j$ are stable (bistable regions), respectively.  The letter `c' after number $i$ denotes the existence of a stable PO that emerged from $\bm{v}_i$. Figure~\ref{fig:allbif} is colored according to the probabilistic calculations described above. Curves are plotted over the colored map, and correspond to transcritical bifurcations (dashed curves), Hopf bifurcations (solid curves), and saddle-node bifurcations of PO (dotted curves) captured in Fig.~\ref{fig:bif_diag}. The division of the parameter space, which is obtained by the two methods independently, into regions of various dynamics yields an interpretative map of the dynamical system.

As shown in Fig.~\ref{fig:allbif}, small values of infection rates lead to the absence of virus populations ultimately for all possible initial conditions; see region (1). Differences in the initial viral population sizes do not affect the prevalence of the more infectious strain, which is accompanied by elimination of the less infectious strain. A general rule that can be read as ``higher infectiousness implies higher prevalence'', as it will be shown further, remains true almost throughout the whole map, including the region of chaos. Nevertheless, despite the similarities of the strains, under the assumptions of the model, the presented map illustrates evident asymmetry.

When infection rates of the generalist and the specialist strains are equal, then, depending on the degree of infectiousness, one or the other strain can persist: for comparatively small values of infection rates, the specialist virus outcompetes the generalist virus. However, for comparable and high infection rates, the generalist prevails. Values of infection rates that lead to the sole prevalence of the generalist strain correspond to large regions, (5) and (5c), on the stability map. Interestingly, the conditions for prevalence of the generalist strain in the system are so favorable that persistence of the generalist affects the stability of other states and leads to bistability. The bistability regions are all found between $T_{57}$ and $H_{5}$ -- the two bifurcations that bring steady stability to the spec-free state. Between these two lines, the spec-free state is simultaneously stable with other states, namely: the virus-free state (1+5), gen-free states (2+5) and (6+5), the coexistence state (3+5), and POs about the gen-free state (6c+5). Among the bistable regions, the (2+5) and (6+5) deserve a special noting. The existence of these regions shows that the generalist or specialist can outcompete each other to become the singular strain surviving in the system for the exact same values of parameters. 
Despite the strong influence of the generalist strain, even the least infectious specialist strains can persist in the bistable regions. However, the specialist with low infection rates can only persist in coexistence with the generalist (see region (3+5)).

\section{Chaos}
\label{sec:chaos}
\begin{figure}
    \subfigure[$\protect \malps = 0.7$]{\includegraphics[width=1\linewidth]{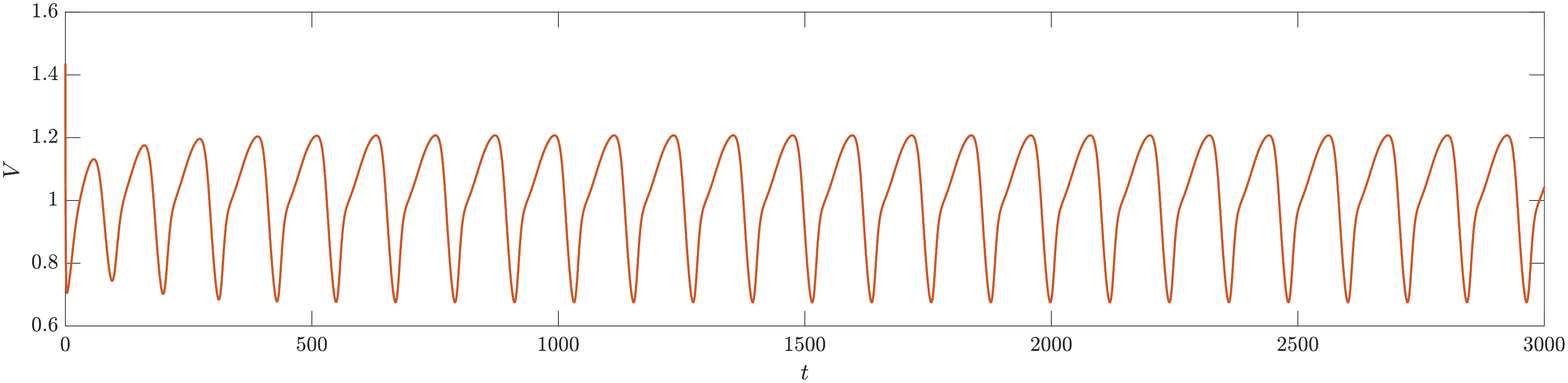}}
    \subfigure[$\protect \malps = 1.2$]{\includegraphics[width=1\linewidth]{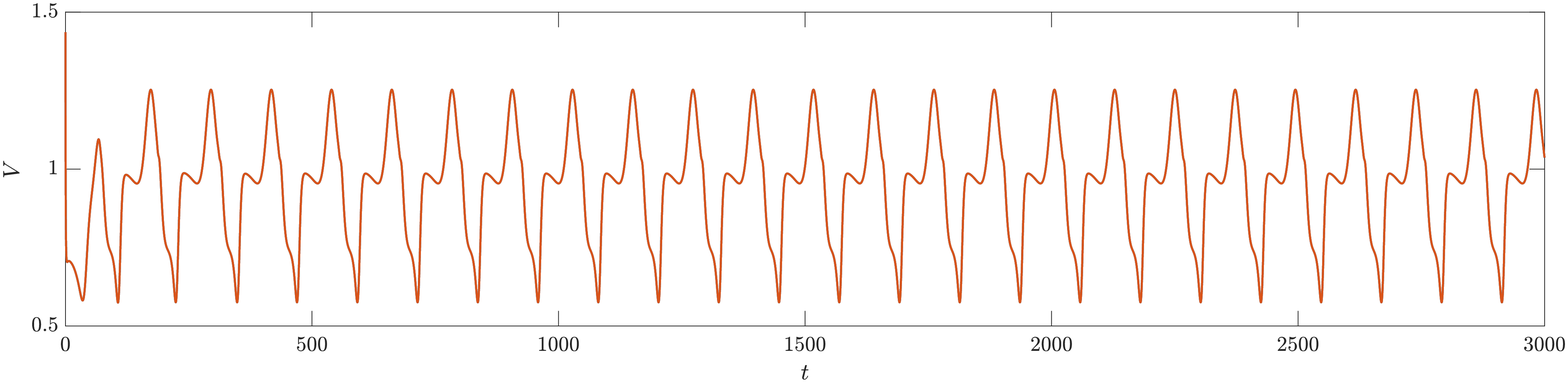}}
    \subfigure[$\protect \malps = 1.6$]{\includegraphics[width=1\linewidth]{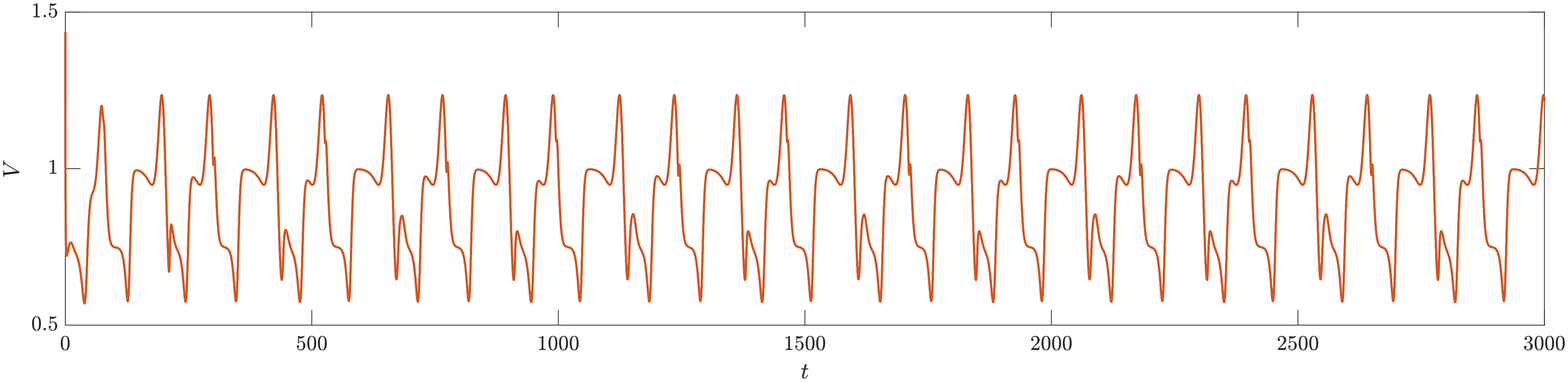}}
    \subfigure[$\protect \malps = 2.0$]{\includegraphics[width=1\linewidth]{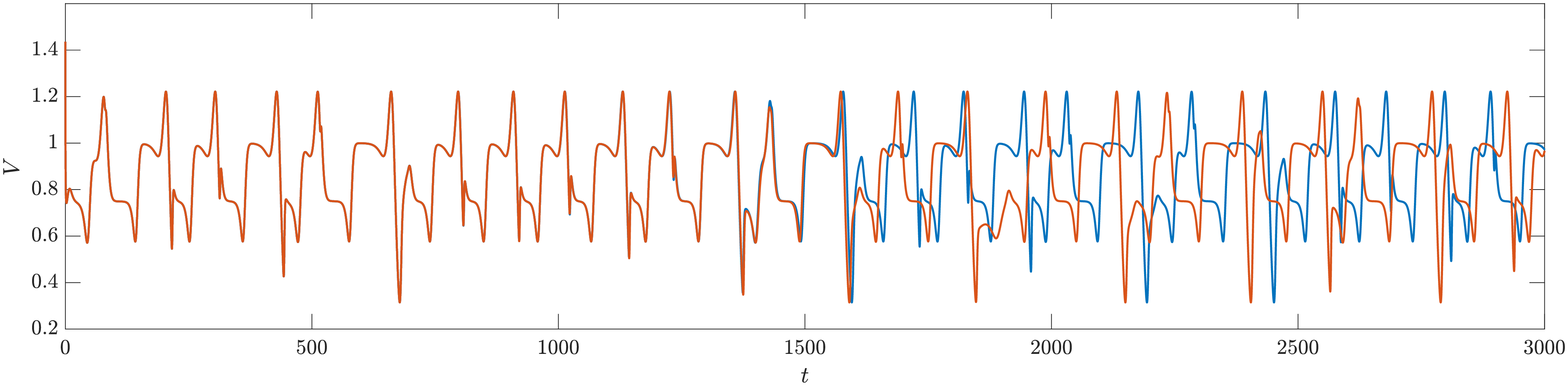}}
	\caption{Time series of system \eqref{odes_twocell_spec_gen_nd} for four different values of \protect \alps. Infection rate $\protect \malpg = 0.5$ and other parameters are fixed as \eqref{params7d}. Initial conditions used for simulations: $(x_1(0) = 1,\, x_2(0) = 0.75,\, \protect \mys(0) = 0,\, \protect  \myg(0) = 0,\, \protect  \mygg(0) = 0,\, \protect  \mzs(0) = 0.5,  \, \protect \mzg(0) = 0.5)$ in red; and $(x_1(0) = 1, \,x_2(0) = 0.75,\, \protect \mys(0) = 0, \,\protect  \myg(0) = 0, \,\protect  \mygg(0) = 0,\, \protect  \mzs(0) = 0.5001,\,  \protect \mzg(0) = 0.5001)$ in blue. The two trajectories from these initial conditions are overlapping in (a), (b), (c). Note the sensitive dependence on initial conditions for $\protect \malps = 2.0$.}
	\label{fig:timeseries_7d}
\end{figure}

In order to understand the dynamics of the system for small \alpg and large \alps where no stable solutions were observed in the one-dimensional bifurcation diagrams, we compute the time series for a fixed value of $\malpg = 0.5$ and four, gradually increasing values of \alps at $0.7$, $1.2$, $1.6$, and $2$.  The results of these computations are shown in Fig.~\ref{fig:timeseries_7d}. Each plot illustrates the time series starting at two, biologically meaningful initial conditions where the second initial condition is a very small perturbation of the first.

At $\malps = 0.7$, the system corresponds to region (6c) of Fig.~\ref{fig:allbif}. These parameters lead to a stable periodic orbit around the gen-free-2 state $\bm v_6$. At $\malps = 1.2$, the periodic orbit shows irregularity which is evident from the emerging peaks, and for $\malps = 1.6$, period doubling becomes apparent. Finally, at $\malps = 2.0$, two major changes are observed.  Firstly, the cycles clearly become aperiodic.  Secondly, the two numerical solutions diverge at the time horizon of $t \approx 1500$. These numerical results thus seem to indicate the presence of chaotic dynamics, which emerge for high values of the infection rate of the specialist and comparably low infection rate of the generalist.  We hypothesize that a cascade of period-doubling bifurcations of periodic orbits is responsible for the onset of chaos. 

To study the chaotic behavior in more detail, we calculate the largest Lyapunov exponents (LLEs, following the method described in Ref.~\cite{parker2012practical}, pages 74--80) as functions of the infection rates.  Positive LLEs are used to 
approximate boundaries that separate chaotic and non-chaotic regions of parameter space.  The red region in Fig.~\ref{fig:allbif} indicates where the LLE is positive with corresponding values shown in the color bar.  Interestingly, the LLEs are only positive in a finite range of infection rates given roughly by $0 < \malpg \leq 0.7$ and $1.5 \leq \malps \leq 3$, which can be considered a chaotic window. 
Outside of this window, the long-time behavior of the system converges to a stable equilibrium point or PO.  For $\malps > 3$, the periodic solution that is born from $\bm{v}_6$ becomes stable.  

\begin{figure}
\centering
{\includegraphics[width=.63\textwidth]{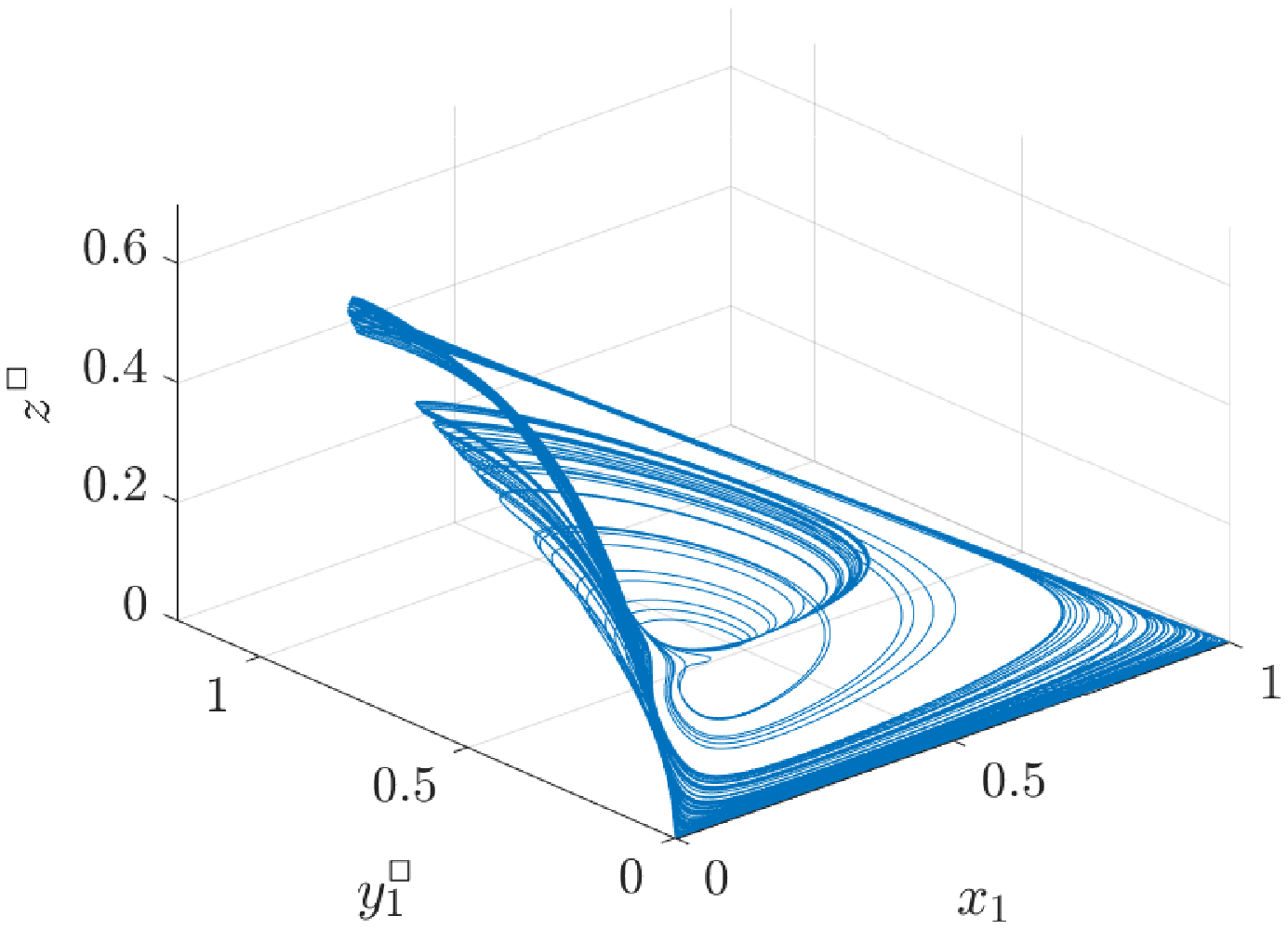}}
{\includegraphics[width=.63\textwidth]{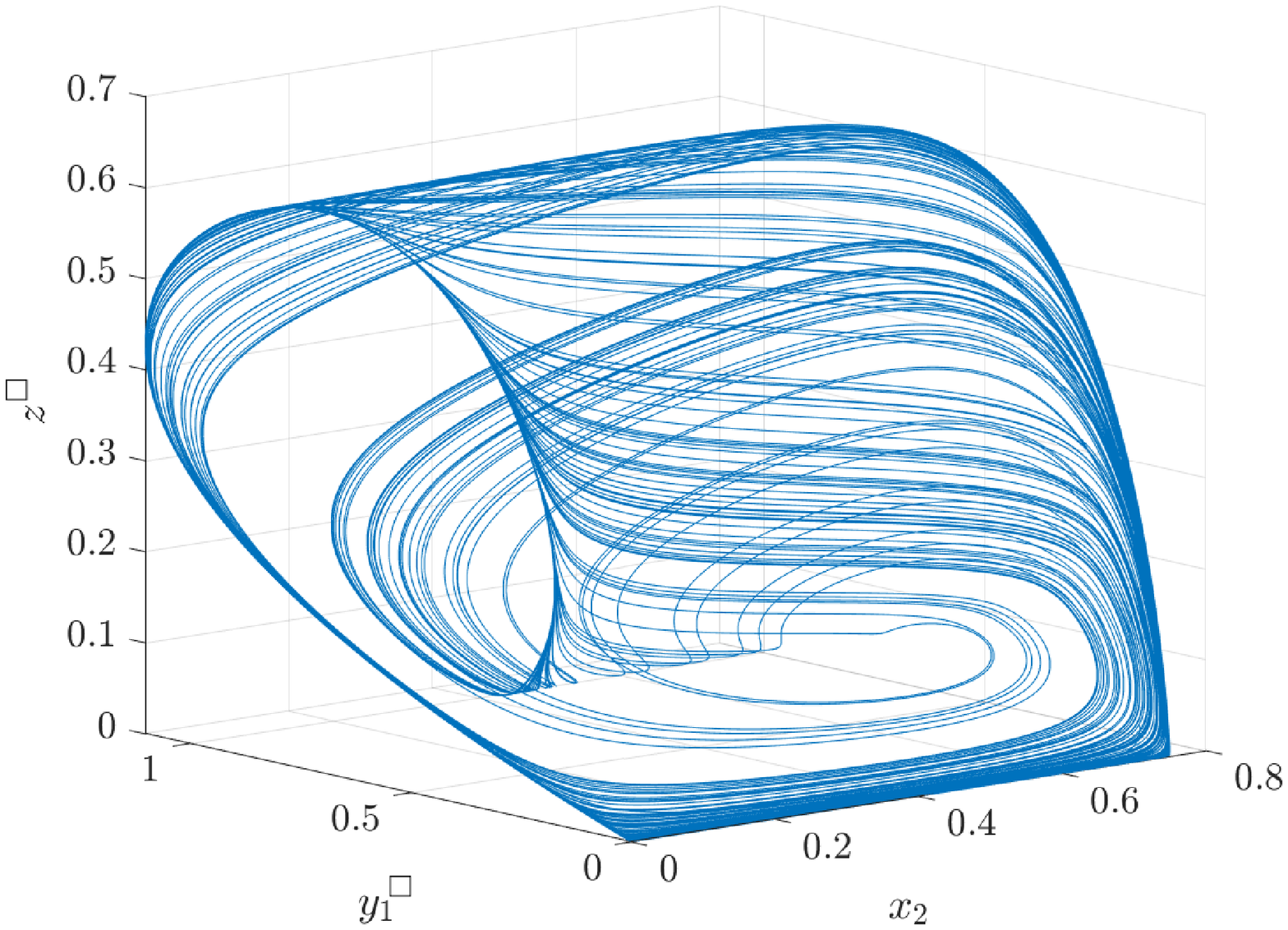}}
{\includegraphics[width=0.6\textwidth]{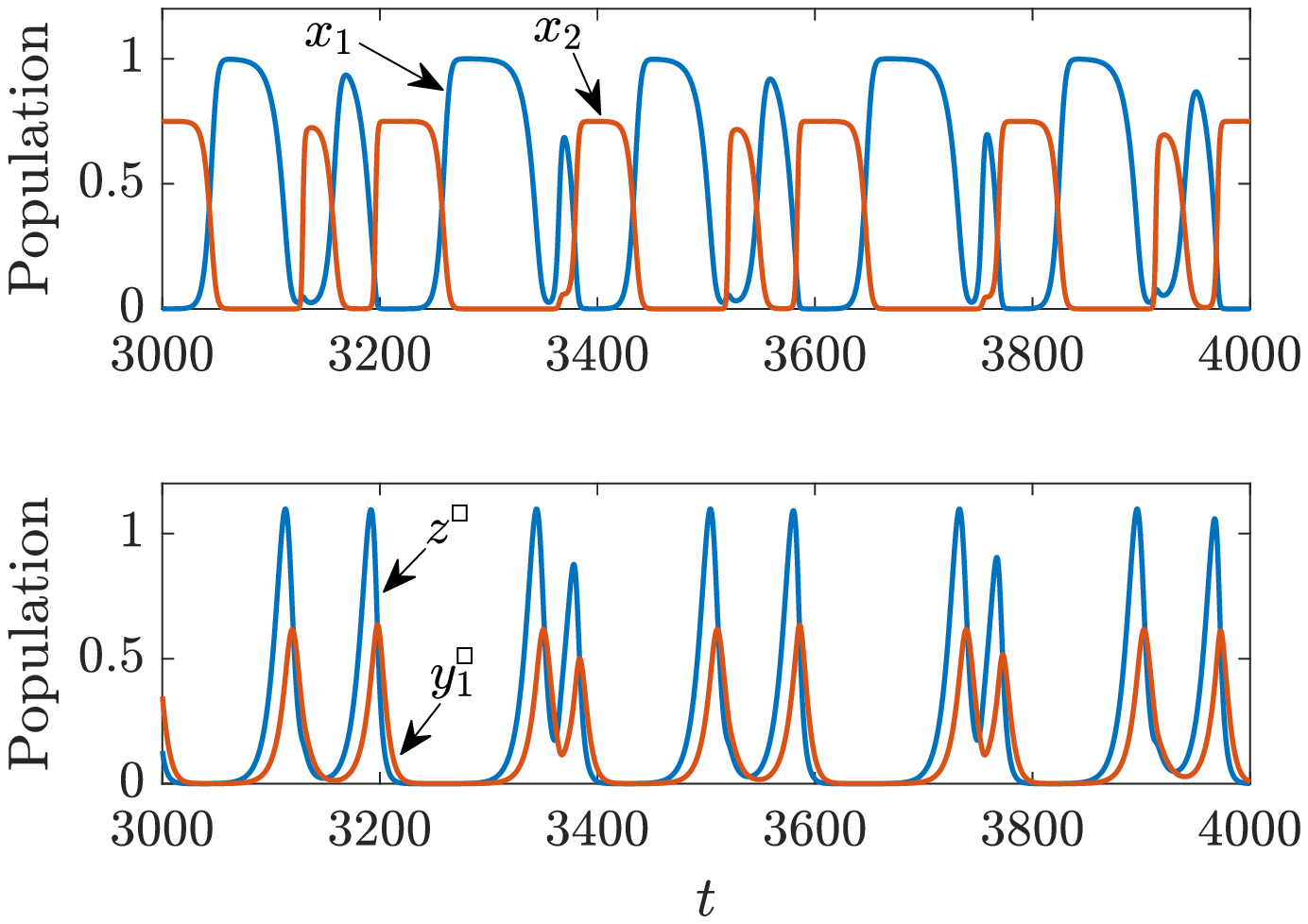}}
	\caption{Three-dimensional projections of the strange chaotic attractor governing the dynamics of Eqs. \eqref{odes_twocell_spec_gen_nd} using parameter values in \eqref{params7d} with $\alpha = 0.5$, $\alpha^{\protect \bbox} = 2.0$.  Below we display time series at the same parameter values showing chaotic  oscillations.	
	}
	\label{fig:strange}
\end{figure}

Projections of the strange chaotic attractor are shown in Fig.~\ref{fig:strange}.  The populations associated with the generalist, \emph{i.e.}\ $\myg$, $\mygg$, and $\mzg$, tend to zero in the chaotic window and thus the dynamics can be captured by a four-dimensional model describing the infection of susceptible cells ($x_1$) that are competing with unsusceptible cells ($x_2$) for the same resource.  The trajectories of these attractors are also plotted as time series in Fig.~\ref{fig:strange} and correspond to chaotic oscillations that can be decomposed into four stages:
\begin{itemize}
\item{The susceptible cells $x_1$ are nearly at their maximum population size while the unsusceptible cells $x_2$, infected cells $\mys$,  as well as the viral populations $\mzs$ are nearly extinct.}
\item{The viral and infected-cell populations experience rapid growth, substantially reducing the population of susceptible cells until they are nearly extinct.  This enables the unsusceptible cells to flourish, but the viral population is quickly diminished as they have exhausted the supply of their hosts.}
\item{The unsusceptible cells are nearly at their maximum population size and the susceptible cells, infected cells, and virus are nearly extinct.}
\item{The remaining susceptible cells rapidly out-compete the unsusceptible cells and the process then repeats.}
\end{itemize}
Due to the chaotic nature of these oscillations, these stages occur at irregular intervals and might only partially develop before transitioning into the next.



\section{Discussion}
\label{sec:conc}
This work provides a comprehensive investigation of an ODE-based model of specialist and generalist viral populations placed in an environment composed of two types of susceptible host cells. The virus can either specialize towards the infection of only one particular type of cell or become a generalist that infects both types (see Fig.~\ref{diag}).  The main goal of this study is to develop a mathematical model  that describes the evolution of such system and  to  determine the dependency of its behavior on the parameters and on the initial conditions focusing on the infection rates of the generalist and specialist strains.  Overall, the system can stabilize at either eight of the nine different equilibria or at periodic orbits around some of them. In addition to that, there are values of parameters that lead to chaos. 

Our study confirms previously known dynamics such as the extinction of both specialist and generalist viral strains when their infection rates are very low. Likewise, the model predicts the existence of a single strain when the infection rate of that particular strain is considerably greater than the infection rate of the other strain; this out-competition is consistent with the results of simpler models that consider competing viral strains in the presence of one type of susceptible cell~\cite{nurtay2019theoretical, fredrickson1981microbial}. 

The trade-off hypothesis assumes that the generalist viral genotype will never dominate and outcompete a population of specialist viruses in their local host~\cite{bedhomme2015emerging, elena2017local} following the adagio that the ``jack-of-all trades is the master of none''.  In sharp contrast with this prediction, experimental evolution of RNA viruses has provided many examples of the evolution of no-cost generalists when the evolving viral populations are equally likely to infect the alternative hosts, see \emph{e.g.} Refs~\cite{turner2000cost, bedhomme2011multihost, deardorff2011west, cooper2001differential}.

Our model provides a theoretical framework to understand these observations.  Although outcompetition of the generalist virus is the most likely outcome of our model (when the infectiousness of the specialist is greater than that of the generalist), we still find regions in parameter space where generalist viruses can coexist with specialists and even outcompete them.  For example, we have observed that the existence of generalist strains whose infection rates lay between $T_{57}$ and $H_{5}$ (exact values are obtainable for every set of parameters) does not depend on the infectiousness of a specialist strain. Counterintuitively, under these conditions, the persistence of the specialist strain is only possible if its infection rate is comparatively \emph{low} compared to the one of the generalist. Interestingly, when \alpg is found between the conditions determined by $\malpg > T_{57}$ and $\malpg < H_{5}$, with \alps being slightly greater than zero, it might result in the coexistence of specialist and generalist strains. From there, a slight increment of the infection rate of \zs, can result in the dominance of the specialist and outcompetition of the generalist. However, for the same parameter values, a change of initial population sizes can result  in the generalist dominating the population and outcompeting the specialist. When \alpg is between $T_{57}$ and $H_{5}$, greater values of infection rate of the specialist, surprisingly, lead only to stronger persistence of the generalist strain. This inverse dependence of survival of the specialist on its infection rate can be explained by the high infectiousness of the generalist. The generalist competes with the specialist for the susceptible cells and, apparently, there might not be enough susceptible cells to maintain persistence of the specialist. Therefore, in the presence of such infectious generalist, high infection rate of the specialist means low overall success of the strain.

One of the main differences between specialist and generalist strains, in terms of system dynamics, is that specialists persist in the chaotic region. The window of chaos is adjacent to periodic persistence of the specialist and periodic coexistence of both strains. Despite the latter, for values of infection rates that force chaotic behavior, the generalist strain is simply not fit enough to persist. To rationalize the chaotic window, note that when the infection rate of the generalist strain is comparatively low, the dynamics of the system are mostly driven by specialist infections. With increasing infectiousness of the specialist, both generalist-free equilibria of the system gain periodic orbits around them. With two fluctuating attractors, the outcome of the system becomes sensitive towards the availability of generalist-susceptible cells, whose value (zero or non-zero) determines the qualitative difference between the two generalist-free equilibria. However, the specialist cannot interact with the generalist-susceptible cells, and an abundance of the latter only depends on presence and infectiousness of the generalist. When the infection rate of the generalist becomes sufficiently high, the system stabilizes with coexisting (and fluctuating) strains.  Right before acquiring the stable periodic orbit around the coexistence state for increasing infection rate of the generalist, amplitudes of fluctuations can be already large due to components related to the specialist string (specialist-susceptible cells, specialist-infected cells, and the specialist strains themselves).

It is expected that a generalist strain flourishes in an environment with multiple host types in comparison to specialist strain that ``prefers'' the environment to which it is most adapted. This would mean that only dynamics from an upper triangle corresponding to values of infection rate $\malpg < \malps$ in the stability diagram (Fig.~\ref{fig:allbif}) would have a chance to be observed in reality. This would, nevertheless, include all the qualitatively different dynamics: (1) virus-free state; (2), (6) generalist-free state; (5) specialist-free state; (6c) periodic generalist-free state; (5c) periodic specialist-free state; (7), (7c) coexistence and periodic coexistence states; and, finally, chaotic behavior. However, as a study on a multi-host disease caused by canine distemper virus has shown~\cite{nikolin2012antagonistic}, there are cases when ``the fitness of generalists should be approximately the same in all environments, including the environment to which a specialist is adapted''. Therefore, the full bifurcation diagram, like the one obtained in this work, can be of a high value when it comes to understanding the dynamics with specialist and generalist infection strategies.



\section*{Acknowledgements}
AN received funding from the ``La Caixa'' Foundation and The Mathematics for Industry Network COST Action [TD1409].  MGH received funding from the European Union's Horizon 2020 research and innovation program under the Marie Sk{\l}odowska-Curie grant [agreement No.~707658]. JS has been funded by a ``Ram\'on y Cajal'' contract [RYC-2017-22243], and by the MINECO grant [MTM2015-71509-C2-1-R] and the Spain's ``Agencia Estatal de Investigaci\'on'' (AEI)  grant [RTI2018-098322-B-I00]. LlA has been supported by the AEI grant [MTM2017-86795-C3-1-P]. SFE support comes from AEI-FEDER grant [BFU2015-65037-P] and Generalitat Valenciana grant [PROMETEU/2019/012]. The research leading to these results has received funding from ``la Caixa'' Foundation, from a MINECO grant awarded to the Barcelona Graduate School of Mathematics (BGSMath) under the ``Mar\'ia de Maeztu'' Program [grant MDM-2014-0445], and from the CERCA Programme of the Generalitat de Catalunya.

\bibliography{short_lib}
\bibliographystyle{unsrt}

\titleformat{\section}{\large\bfseries}{\appendixname~\thesection .}{0.5em}{}
\begin{appendices}
\section{Computation of Equilibria in Table~\ref{TheCPTable} }\label{app:CP}

\paragraph{{\bfseries Row \boldmath{$v_0$}.}}
In the conditions given in Column~3 of Table~\ref{TheCPTable}
for $v_0$, we get $z^* = 0$ from the third equation of~\eqref{TheSystem} because $\zeta \neq 0.$
\paragraph{Remark.} {\boldmath{$x_1 = 0 \Longrightarrow \mzs = 0$}} by using the fourth equation of~\eqref{TheSystem}
together with the fact that $\zeta^{\bbox} \neq 0.$
This determines the seventh column of Table~\ref{TheCPTable} for Equilibria~$v_0,\ v_4,$ and $v_5.$
\paragraph{{\bfseries Row  \boldmath{$v_1$}.}}  It follows trivially from the first equation of~\eqref{TheSystem}.
\paragraph{{\bfseries Row  \boldmath{$v_2$}.}}
Since $\mzs \neq 0,$ the value of $x_1^*$ is obtained from the fourth equation of~\eqref{TheSystem} after plugging $z=0$ to it. Then, $\mzs{}^*$ can be computed from the first
equation of~\eqref{TheSystem} under the conditions of this row 
(and after plugging the value of $x_1^*$ into it).
\paragraph{{\bfseries Row  \boldmath{$v_3$}.}} 
By plugging the conditions of this row into the third and the first equations of~\eqref{TheSystem},
we obtain $x_1^*$ and $z^*,$ respectively.
Then, by plugging the already computed value of $z^*$ into the fourth equation of~\eqref{TheSystem},
we get
$1-x_1 = \mzs\left(\tfrac{\alpha\zeta^{\bbox}}{C x_1} + \malps + \tfrac{D \alpha}{C}\right)$,
which gives $z^{\bbox}{}^*$ in this case.
\paragraph{{\bfseries Row  \boldmath{$v_4$}.}} 
The component $x_2^*$ of the critical point can be obtained from the second equation of~\eqref{TheSystem}.
\paragraph{{\bfseries Row \boldmath{$v_5$}.}}    
Since $x_2,z \neq 0,$ the third and the second equations of~\eqref{TheSystem} give  $x_2^*$ and $z^*,$ respectively.
\paragraph{{\bfseries Row \boldmath{$v_6$}.}} 
Observe that {\boldmath{$z = 0 \Longrightarrow \mzs \neq 0$}}.
Indeed, when $z = \mzs = 0$ the third and the 
fourth equations of~\eqref{TheSystem}
are verified automatically, and the first and 
the second equations read $x_1 = 1 - x_2$ and
$
0 = \beta_1 - \beta_2(x_1+x_2) = \beta_1 - \beta_2 \neq 0.
$
Thus, (since $z=0$ and $\mzs \neq 0$) the fourth equation 
gives $x_1^*.$ Next, we obtain $x_2^*$ from the second equation,
after plugging $z=0$ and the already known value of $x_1^*$ into it.
Finally, again from the second equation (or from the two components already computed), 
we obtain $x_1 + x_2 = \tfrac{\beta_1}{\beta_2}.$ By plugging this value and $z=0$
into the first equation, we get $\mzs{}^*.$
\paragraph{{\bfseries Rows \boldmath{$v_7$}} and {\boldmath{$v_8$}.}} 
From the second and the third equations of~\eqref{TheSystem} we obtain the values of $z^*$ and $x_2^*,$
respectively. Then, the value of $\mzs{}^*$ can be obtained from the first
equation after plugging the already computed expression of $z^*$ into it.
Observe that this gives the values of $x_2^*,\ z^*,$ and $\mzs{}^*$ in terms
of the (for now) unknown value of $x_1^*.$ To determine the latter, we plug the three
expressions for $x_2^*,\ z^*,$ and $\mzs{}^*$ into the fourth equation to get
\[ \frac{\phi_2 x_1^2 + \phi_1 x_1 + \phi_0}{B\alpha\malps} = 0, \]
where
\begin{align*}
\phi_2 & = (A-B)\Bigl(\bigl(D\alpha + C\malps\bigr)\beta_2 - D\alpha\Bigr), \\
\phi_1 & = (A-B)\alpha\zeta^{\bbox}\bigl(\beta_2	- 1\bigr) + 
           \bigl(D\alpha + C\malps\bigr)\bigl(B\beta_1 - \zeta\beta_2\bigr) + D \alpha(\zeta-B),\text{ and}\\
\phi_0 & = \alpha\zeta^{\bbox}\Bigl(\bigl(1-\beta_2\bigr)\zeta + B\bigl(\beta_1 - 1\bigr)\Bigr).
\end{align*}
Since $B\alpha\malps \ne 0$, the equation above gives two critical values of $x_1$, denoted by 
$x_1^{\pm}$, which are the two solutions of the quadratic equation 
$\phi_2 x_1^2 + \phi_1 x_1 + \phi_0,$
whenever they are real and non-negative.

\end{appendices}

\end{document}